\documentclass[12pt,a4paper]{article}

\usepackage[T1]{fontenc}
\usepackage{lmodern}
\usepackage{amsmath}
\usepackage{amsthm}
\usepackage{url}
\usepackage{makecell}
\usepackage{setspace}
\usepackage{graphicx}

\usepackage[a4paper,left=3cm,right=3cm,top=3cm,bottom=3cm]{geometry}

\newcommand{\eqb}{\begin{equation}}
\newcommand{\eqe}{\end{equation}}
\newcommand{\al}{\alpha}

\title{From geometry to sustainability: Optimal shapes of hip roof houses}

\author{
Ewa Rokita-Magdziarz\thanks{
Rokita-Project Architectural Office, Forsycjowa 7, Wroc{\l}aw, 51-253, Poland.\\
Corresponding author: \texttt{biuro@rokita-projekt.pl}
}
\and
Barbara Gronostajska\thanks{
Department of Architectural and Construction Design, Faculty of Architecture,
Wroc{\l}aw University of Science and Technology, Wyspia\'nskiego 27, Wroc{\l}aw, 50-370, Poland.
}
\and
Marcin Magdziarz\thanks{
Department of Applied Mathematics, Faculty of Pure and Applied Mathematics,
Wroc{\l}aw University of Science and Technology, Wyspia\'nskiego 27, Wroc{\l}aw, 50-370, Poland.
}
}

\date{} 

\begin{document}

\maketitle

\begin{abstract}
\begin{spacing}{0.9}
 In this paper, we develop a rigorous mathematical framework for the optimization of hip roof house geometry, with the primary goal of minimizing the external surface of the building envelope for a given set of design constraints. Five optimization scenarios are systematically analyzed: fixed volume, fixed footprint ratio, fixed slenderness ratio, fixed floor area, and constrained height. For each case, explicit formulas for the optimal dimensions are derived, offering architects and engineers practical guidelines for improving material efficiency, reducing construction costs, and enhancing energy performance. To illustrate the practical relevance of the theoretical results, case studies of real-world hip roof houses are presented, revealing both inefficiencies in common practice and near-optimal examples. Furthermore, a freely available software application has been developed to support designers in applying the optimization methods directly to architectural projects. The findings confirm that square-based footprints combined with balanced slenderness ratios yield the most efficient forms, while deviations toward elongated or flattened proportions significantly increase energy and material demands. This work demonstrates how mathematical modeling and architectural design can be integrated to support sustainable architecture, providing both theoretical insight and practical tools for shaping energy-efficient, cost-effective, and aesthetically coherent residential buildings.
\end{spacing}
\end{abstract}


\maketitle
\newpage
\section{Introduction}

The geometry of a building is one of the most decisive factors influencing its structural performance, material efficiency, and energy demand. Among the variety of residential typologies, the hip roof house represents a particularly important case, both in traditional and contemporary contexts. Its geometry, defined by four sloping planes converging on a ridge or peak, not only provides aesthetic harmony but also ensures structural resilience under wind loads, efficient rainwater drainage, and balanced load distribution on walls and foundations. These inherent advantages have made hip roof houses a popular and enduring choice in different regions, and surveys consistently rank them among the most common roof forms worldwide \cite{hip_roof_www}. At the same time, however, the design of hip roof houses is often guided more by convention, stylistic preference, or regulatory requirements than by systematic optimization. As a result, many built examples deviate significantly from their theoretical optima, leading to avoidable increases in material use, construction cost, and energy consumption.

The connection between building form and energy efficiency has been highlighted by a wide range of studies. Depecker et al. \cite{depecker2001design} established the influence of building shape on heating demand, showing that more compact envelopes minimize thermal losses. Hegger et al. \cite{hegger2012energy} stressed the same principle in their design manual on sustainable architecture, reinforcing the role of geometry in reducing operational energy. Similarly, Ramesh et al. \cite{Ramesh2010} demonstrated in a comprehensive life-cycle energy analysis that form-related parameters can account for a substantial fraction of total building energy use. Numerous studies have attempted to model and quantify this relationship, such as Catalina et al. \cite{Catalina2011}, Ourghi et al. \cite{Ourghi2007}, Schlueter and Thesseling \cite{Schlueter2009}, and Jedrzejuk \cite{Jedrzejuk2000, Jedrzejuk2002}, all of which confirmed the strong correlation between compactness and performance. Steemers \cite{Steemers2003} further argued that architectural design at the morphological level is one of the most effective means for reducing energy demand. More recent contributions by D’Amico and Pomponi \cite{d2019compactness} introduced a dimensionless compactness measure that allows scale-independent evaluation of building geometries, overcoming the limitations of traditional surface-to-volume ratios. Their work demonstrates that geometry is as decisive for sustainability as material choice or systems design.

Urban and typological studies reinforce these findings at larger scales. Steadman’s analyses of design evolution and building types \cite{Steadman2000,Steadman2014} highlight the enduring influence of form on efficiency, while Julia et al. \cite{Julia2017} and Hargreaves \cite{Hargreaves2017} incorporated typological and density parameters into urban energy modeling. Angel’s \cite{Angel2012} global perspective on urban expansion, together with UN reports \cite{UN2014} and policy frameworks such as the London Strategy \cite{LondonStrategy2018} and California’s energy codes \cite{California2017}, emphasize that residential design decisions have long-lasting implications for sustainable urban development. The hip roof house, as a representative typology, thus plays a role not only in architectural practice but also in shaping the energy and carbon trajectories of cities.

The theoretical background of architectural form generation provides another foundation for this research. Early contributions  \cite{Martin1972,March1976,Stiny1980} pioneered the study of shape grammars, showing that design can be formalized and reproduced through logical rules. Mitchell \cite{Mitchell1990} deepened this connection by embedding formal logic into computational design processes, while Knight \cite{Knight1994} expanded the methodology into a robust system for shape exploration. These studies, synthesized by Oxman \cite{oxman2006theory}, underpin much of contemporary computational design thinking, which sees architecture as a fusion of generative rules and performance-driven criteria. Within this tradition, geometry is no longer incidental but becomes an explicit design parameter that can be optimized alongside other requirements. 

Optimization itself has increasingly become a cornerstone of sustainable architecture established convex optimization as a rigorous framework for solving design problems under constraints, offering powerful tools for deriving global optima in form-finding \cite{boyd2004convex}. Okeil \cite{Okeil2010}, Caruso et al. \cite{Caruso2013}, Jin and Jeong \cite{Jin2017}, and Vartholomaios \cite{Vartholomaios2017} applied such methods in architectural design, ranging from parametric explorations to performance-based optimization. Hachem \cite{Hachem2012,Hachem2016} extended the approach to neighborhood scales, where residential morphologies are optimized to maximize solar potential and energy efficiency. These works confirm that optimization methods are not abstract mathematical exercises but practical tools for achieving sustainability in the built environment.

In light of these contributions, the present paper proposes a rigorous mathematical framework for the optimization of hip roof houses. By deriving closed-form formulas for optimal geometry under diverse constraints—fixed volume, footprint aspect ratio, slenderness ratio, floor area, and limited height—we aim to demonstrate how compactness theory, optimization methods, and typological research can be integrated into practice. Case studies of existing hip roof houses are used to validate the results, while a dedicated software application enables practitioners and students to directly apply the methods. In doing so, this research bridges the gap between theoretical models and architectural practice, contributing to both sustainable housing design and the broader discourse on the role of geometry in the built environment.

In this paper, we develop and apply a rigorous optimization framework for the geometry of hip roof houses. Section~2 presents the mathematical formulation of the optimization problem and derives explicit solutions for five design scenarios, namely: (i) fixed volume, (ii) fixed footprint ratio, (iii) fixed slenderness ratio, (iv) fixed floor area, and (v) constrained height within a predefined range. Each case is analyzed in detail, with closed-form expressions provided for the optimal dimensions and minimal external surface. In Section~3, the theoretical results are tested against real-world examples through three case studies of existing hip roof houses, where the derived formulas are applied to evaluate compactness, identify inefficiencies, and quantify potential savings. Section~4 introduces a dedicated software application, developed to facilitate the practical use of the optimization methods in architectural and engineering design workflows, making the theoretical findings directly accessible to practitioners and students. Finally, Section~5 summarizes the results, highlights their implications for sustainable housing design, and outlines possible directions for future research in extending the optimization approach to more complex roof typologies and integrated energy performance assessments.

\section{Optimization framework for hip roof geometry}

Hip roof houses are one of the most popular choices among investors, now representing the second most common style of roof in the U.S. (after a gable roof) \cite{hip_roof_www}.
Hip roof houses combine practical and aesthetic qualities that contribute to their popularity across both traditional and modern architecture. With slopes on all four sides, they offer strong resistance to high winds, which enhances their stability and durability compared to many other roof forms. This geometry also allows rain, snow, and ice to drain efficiently, minimizing the risk of leaks or structural damage and extending the building’s lifespan. Their compact shape supports effective insulation and ventilation, which, when paired with quality materials, can lead to lower energy consumption for heating and cooling. In addition to these functional benefits, hip roofs are highly regarded for their symmetrical, refined appearance that adapts well to a wide range of architectural styles. They also provide opportunities for additional attic or loft space and can be combined with dormers or extensions to enrich the overall design. Owing to their resilience, energy efficiency, and visual appeal, hip roof houses are frequently associated with enhanced property value and strong investor interest.

\begin{figure}[t]
\centering
\includegraphics[width=12cm]{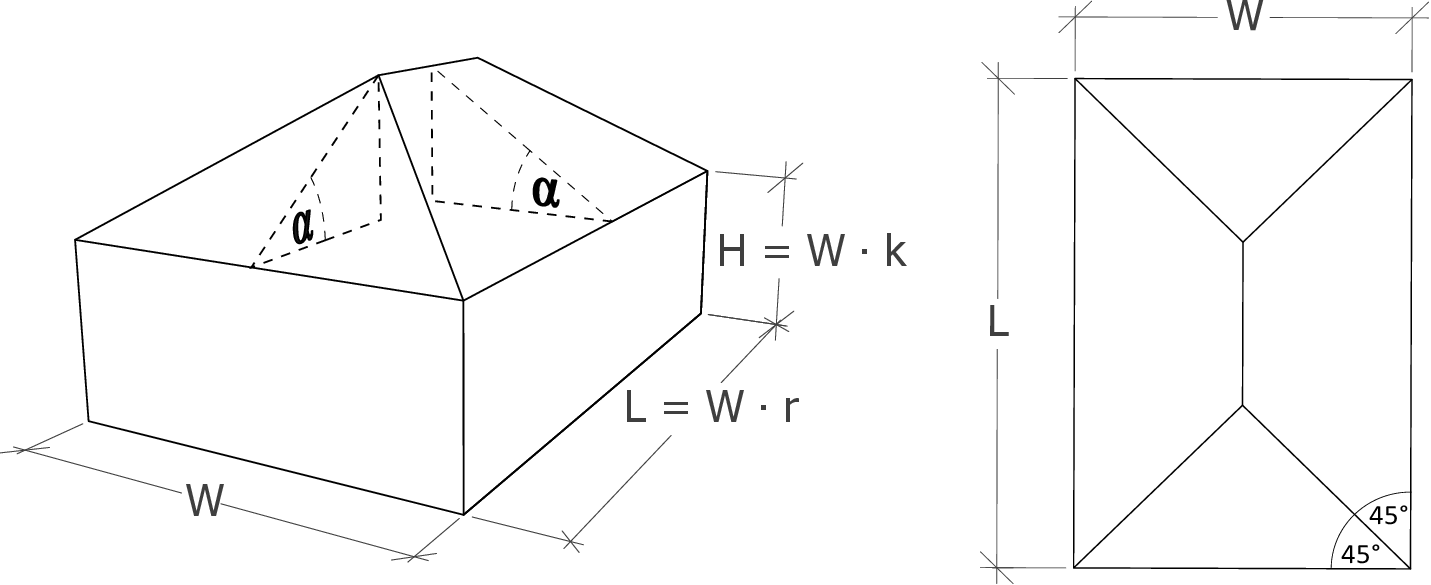}
\caption{Typical shape of the hip roof house with notation used in the paper. On the left side of the drawing we have a perspective view, on the right a view from above. }
\label{koperta}
\end{figure}

A standard form of a hip roof house is illustrated in Fig. \ref{koperta}. The roof consists of two larger, trapezoidal-shaped faces and two smaller, triangular-shaped faces at either end. Assuming equal slope of all roof surfaces, the roof of a building based on a square plan is composed of four identical triangular faces that meet at a single peak, forming a pyramid shape. All four roof planes slope downward at the same angle $\al$ toward the walls.

The geometry of the hip roof house is defined by the following parameters: $W$ — width, $L$ — length, $H$ — height,  $\alpha$ — roof slope angle.
Additionally, we define two aspect ratios: 
\begin{itemize}
    \item $r = \frac{L}{W}$ — footprint aspect ratio;
    \item $k = \frac{H}{W}$ — slenderness aspect ratio.
\end{itemize}
The roof slope angle $\alpha$ is assumed to be constant and its range typically specified in the Local Special Development Plan within a permissible range. The volume of the habitable area of the house is given by:
\begin{equation}
\label{objetosc}
V = W L H.
\end{equation}
For simplicity, it is assumed that the attic space is non-habitable.

The external surface area of the hip roof house, comprising both walls and roof planes (i.e., the building envelope), is expressed as:
\eqb
\label{powierzchnia}
S=2WH+2LH+2L\frac{\frac{1}{2}W}{\cos(\alpha)}=
2WH+2LH+\frac{LW}{\cos(\alpha)}.
\eqe
The first two components represent the area of the walls. The last component is equal to the area of the roof planes.
It should be noted that the area of the building in contact with the ground is not included in $S$.

The primary objective of this section is to determine the optimal shape of the hip roof house, which minimizes the external surface $S$. We will consider five different scenarios.
\subsubsection{First scenario: fixed volume of the house}
\label{Scenario1}

In this scenario we assume that the volume $V$ is fixed (say,  determined at the investor's request). The angle $\al$ is also fixed (specified in the local development plan).
The objective is to determine the parameters $W$, $L$, and $H$ that minimize the external surface of the building $S$. 
To this end, let us first express the surface $S$ in \eqref{powierzchnia} as a function of two ratios, $r$ and $k$. Recall that $L = W r$, $H = W k$, and $V = W L H$, which implies
\begin{equation}
    W = \left( \frac{V}{r k} \right)^{1/3}.
\end{equation}
It follows that
\begin{align}
S=S(r,k) &=2W^2 k + 2 W^2 rk +  \frac{W^2 r}{\cos(\al)}
\nonumber \\
&= W^2 \left(2k + 2rk + \frac{r}{\cos(\al)}\right) \nonumber \\
&= V^{2/3}\left( \frac{2k + 2rk + \frac{r}{\cos(\al)}}{(rk)^{2/3}}\right).
\label{powierzchnia2}
\end{align}

Since $V$ is fixed, to find the minimum of the external surface $S(r,k)$ it is enough to find the minimum of the function in brackets in \eqref{powierzchnia2}. We denote this function by $\gamma(r,k)$
\eqb
\label{gamma_f}
\gamma(r,k)=\frac{2k + 2rk + \frac{r}{\cos(\al)}}{(rk)^{2/3}}.
\eqe
To find possible extrema points of the function $\gamma(r,k)$, let us calculate the partial derivatives $\frac{\partial \gamma(r,k)}{\partial r}$ and $\frac{\partial \gamma(r,k)}{\partial k}$. We get
\[
\frac{\partial \gamma(r,k)}{\partial r}=
\frac{\frac{kr}{\cos(\al)} +2k^2 r -4k^2}{3(rk)^{5/3}},
\]
\[
\frac{\partial \gamma(r,k)}{\partial k}=
\frac{2kr^2 +2kr-\frac{2r^2}{\cos(\al)}}
       {3(rk)^{5/3}}.
\]
To determine the minimum of the function 
$\gamma(r,k)$, we solve the following system of equations
\begin{equation}
\label{system}
\left\{
\begin{aligned}
&\frac{\partial \gamma(r,k)}{\partial r} = 0, \\[6pt]
&\frac{\partial \gamma(r,k)}{\partial k} = 0.
\end{aligned}
\right.
\end{equation}
After performing the necessary calculations, we obtain that for $\alpha \in (0, \pi/2)$, the positive solution of \eqref{system} that minimizes the function $\gamma(r,k)$ is given by
\eqb
\label{r_min}
r_{\min} =1,
\eqe
\eqb
\label{k_min}
k_{\min} =\frac{1}{2\cos(\al)}.
\eqe
Since the determinant of the Hessian matrix and the second derivative $\frac{\partial^2 \gamma(r,k)}{\partial r^2}$ are positive at $r_{\min}$ and $k_{\min}$, these points determine the minimum of the function $\gamma(r,s)$.
The resulting minimal external surface of the hip roof house equals
\eqb
\label{S_min}
S_{\min}=V^{2/3}\gamma(r_{\min},k_{\min})=\frac{3\cdot (2V)^{2/3}}{(\cos(\al))^{1/3}}.
\eqe
The other optimal parameters of the hip roof house are as follows
\eqb
\label{W_min}
W_{\min}=\left( \frac{V}{r_{\min}\cdot k_{\min}} \right)^{1/3}=
(2V\cos(\al))^{1/3},
\eqe
\eqb
\label{L_min}
L_{\min}=W_{\min}\cdot r_{\min}=(2V\cos(\al))^{1/3},
\eqe
\eqb
\label{H_min}
H_{\min}=W_{\min}\cdot k_{\min}=\frac{(2V\cos(\al))^{1/3}}{2\cos(\al)}=\left(\frac{V}{4\cos^2(\al)}\right)^{1/3}.
\eqe
Note that, as expected, the optimal shape of the hip roof house is achieved when the base of the house is a square ($r=1$), since then we obtain the most compact form.

The expressions in \eqref{S_min}--\eqref{H_min} provide a practical framework for designing a hip roof house with the minimal external surface $S_{\min}$ for a specified volume $V$. By applying these formulas, designers can optimize the building’s geometry to reduce material usage and construction costs. Moreover, minimizing the external surface contributes to improved energy efficiency by lowering heat loss through the building envelope, thereby enhancing the building’s overall sustainability and reducing long-term operating expenses.

Let us analyze the following example: take $V=400\; m^3$ and $\al=\pi/6=30^\circ$. Then, from equations \eqref{r_min}--\eqref{H_min} we get the following optimal parameters of the house:
$r_{\min}=1$, $k_{\min}=0,58$, $S_{\min}=271,23\;m^2 $, $W_{\min}=8,85\;m$, $L_{\min}=8,85\;m$, $H_{\min}=5,11\;m$.

It is important to emphasize that, from a practical standpoint, the above obtained parameters are  applicable to real-world house construction. They not only offer feasible design solutions but also ensure structural efficiency, cost-effectiveness, and enhanced energy performance. As such, these results can serve as a reliable foundation for architects and engineers seeking to implement optimized hip roof house designs in practice.

\begin{figure}[t]
\centering
\includegraphics[width=14cm]{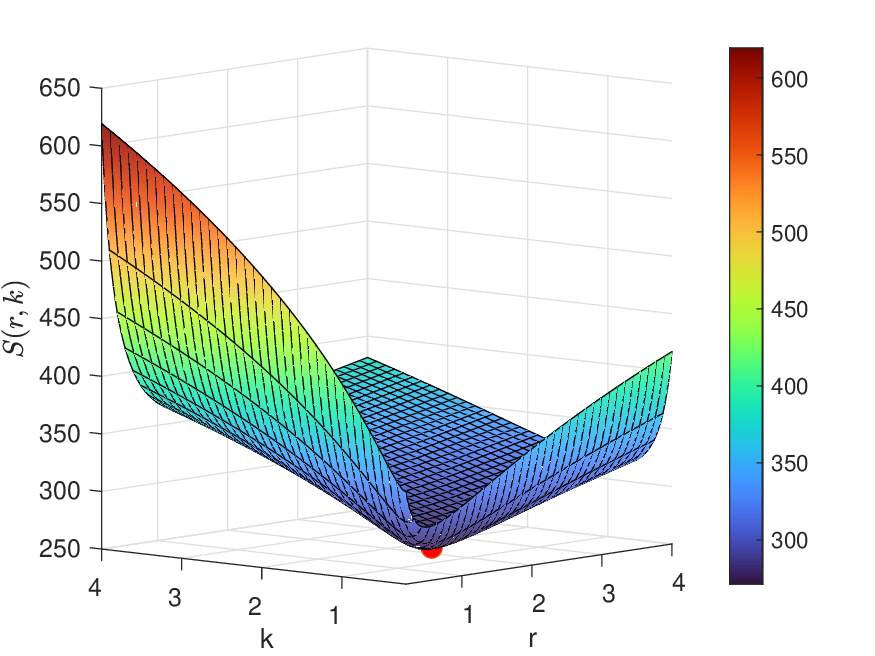}
\caption{Graph of the external surface $S(r,k)$ depending on the footprint aspect ratio $r$ and slenderness aspect ratio $k$, with parameters $V = 400 \; \text{m}^3$ and $\alpha = \pi/6 = 30^\circ$. The red dot indicates the location of the minimal surface. }
\label{S_example}
\end{figure}

Fig. \ref{S_example} illustrates the external surface $S(r,k)$ as a function of the footprint aspect ratio $r$ and the slenderness aspect ratio $k$, for the parameters $V = 400 \; \text{m}^3$ and $\alpha = \pi/6 = 30^\circ$. This graphical representation clearly demonstrates how the surface varies with changes in $r$ and $k$. Clearly, the function has a global minimum. The point corresponding to the minimal surface is highlighted by a red dot, indicating the optimal configuration under the given conditions.

\begin{figure}[t]
\centering
\includegraphics[width=14cm]{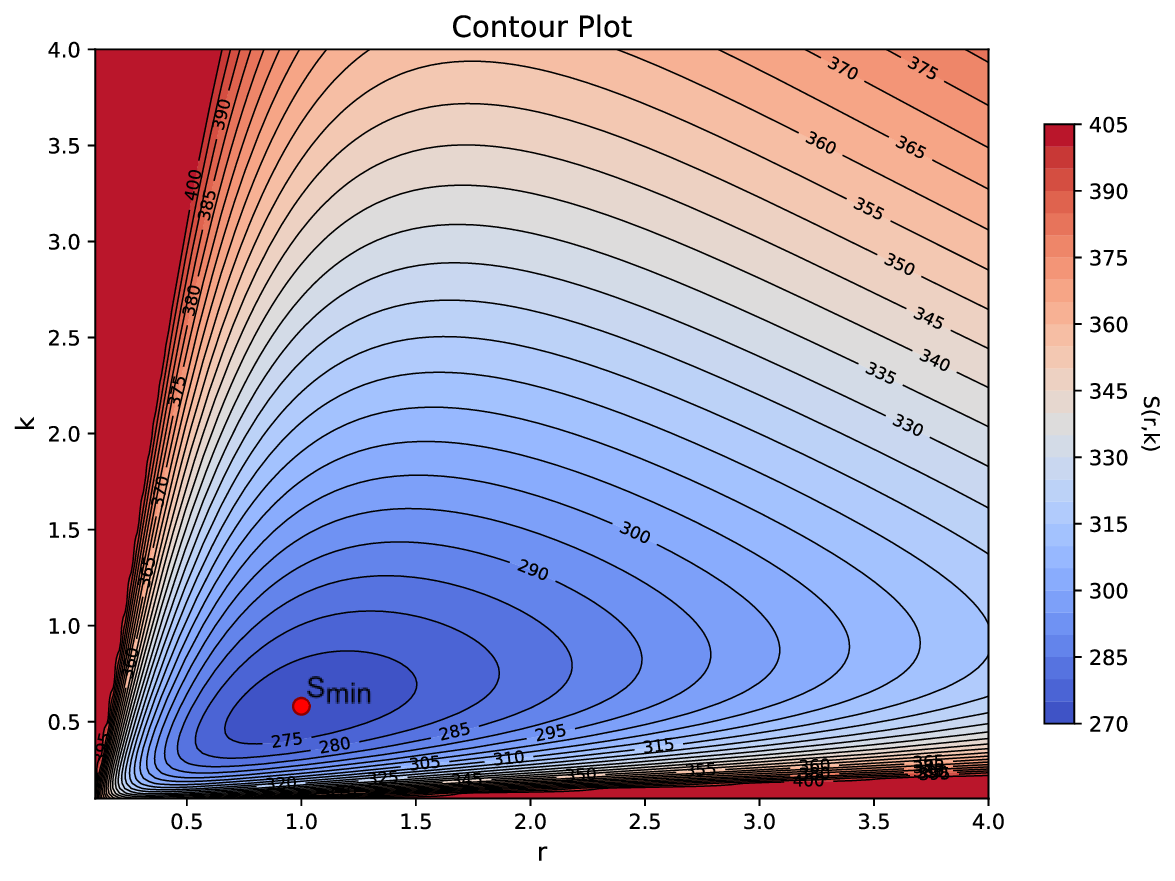}
\caption{Contour plot of the external surface $S$ as a function of  footprint aspect ratio $r$ and slenderness aspect ratio $k$. Here $\al=\pi/6=30^\circ$ and $V=400
\; m^2$. Red dot is the point corresponding to the global minimum equal to $S_{\min}=271,23\;m^2$, attained for $r=1$ and $k=0,58$. }
\label{poziomice}
\end{figure}

Fig. \ref{poziomice} presents the external surface $S$ in the form of a contour plot, using the same parameters as in Figure \ref{S_example}.  
The numerical labels on the contour lines represent the corresponding values of $S$ depending on $r$ and $k$, thereby indicating the degree of deviation from the optimal external surface.  
Each contour line thus identifies the set of $(r,k)$ pairs that yield the same value of $S$.  
This visual representation provides a valuable tool in the preliminary design phase of a hip roof house, enabling designers to readily identify parameter combinations that achieve a specified envelope of the building.  
Such analysis facilitates the selection of an appropriate house shape that not only harmonizes with the constraints of the building plot but also maintains structural efficiency and compactness close to the theoretical optimum.  
Consequently, this approach supports both practical decision-making and optimization in architectural planning.

It is valuable to examine how the optimal parameters $W_{\min}$, $L_{\min}$, and $H_{\min}$ vary with changes in the volume $V$ and the roof slope angle $\alpha$.  
From equations \eqref{W_min}--\eqref{H_min}, it is evident that each expression contains a factor of $V^{1/3}$, which explicitly characterizes the scaling of the optimal parameters with respect to the volume $V$. The parameters increase with increasing $V$, as expected..

\begin{figure}[t]
\centering
\includegraphics[width=10cm]{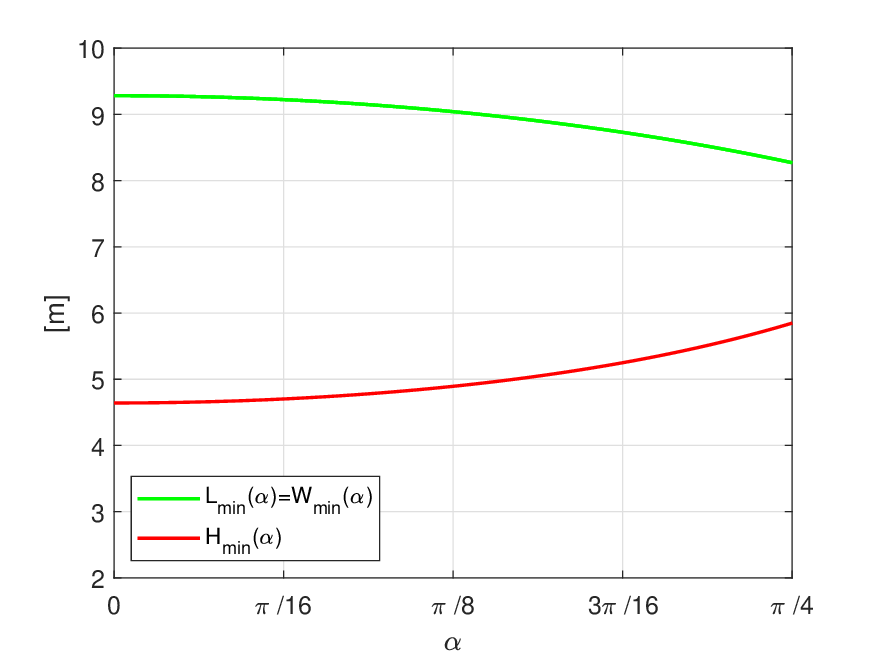}
\caption{Plots of the optimal parameters $W_{\min}$, $L_{\min}$ and $H_{\min}$ as functions of $\al$. Here $V=400 \; m^3$.  }
\label{alpha_dependence}
\end{figure}

A more nuanced behavior emerges when considering the dependence on $\alpha$.  
Fig. \ref{alpha_dependence} presents the behaviour of $W_{\min}$, $L_{\min}$, and $H_{\min}$ as functions of $\alpha$.  
Notably, $W_{\min}$ and $L_{\min}$ exhibit a decreasing trend as $\alpha$ increases, while $H_{\min}$ grows with $\alpha$.  
This analysis provides deeper insight into the interplay between geometric optimization and roof slope, which is of practical significance for architectural design and energy efficiency considerations.

\paragraph{Compactness measure of the building.}
\ \\ 
A widely used approach to characterize the relationship between a building's external surface $S$ and its volume $V$ is the surface-to-volume ratio, $S/V$ \cite{depecker2001design,hegger2012energy}.  
Lower values of $S/V$ are generally preferred, as they indicate greater compactness.  
However, this ratio has an important limitation: it is not scale-independent. That is, two geometrically identical shapes (e.g., cubes) of different sizes will produce different $S/V$ values.  
To address this limitation, a novel dimensionless compactness metric has recently been proposed, designed to rigorously evaluate building compactness based solely on form, thereby removing dependence on absolute size or volume \cite{d2019compactness}.  
This compactness measure is defined as \cite{d2019compactness}:
\begin{equation}
\label{measure_def}
\frac{S}{S_{\min}},
\end{equation}
where $S$ denotes the external surface area of the building, and $S_{\min}$ represents the minimum external surface area required to enclose a given volume $V$.  
This newly introduced metric offers significant value in the early design stage by providing a quantitative assessment of how closely a building's form approaches optimal compactness, thereby enabling informed decisions and guiding the exploration of alternative geometries that may achieve improved efficiency.

\begin{figure}[t]
\centering
\includegraphics[width=14cm]{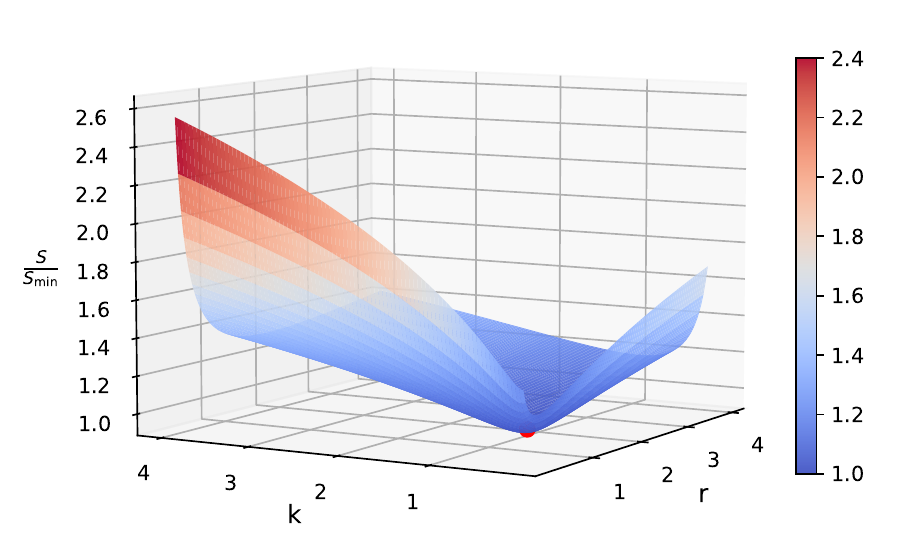}
\caption{Graph of the ratio $S/S_{\min}$ depending on the footprint aspect ratio $r$ and slenderness aspect ratio $k$, with parameters $V = 400 \; \text{m}^3$ and $\alpha = \pi/6 = 30^\circ$. The red dot indicates the location of the minimum equal to 1. }
\label{Compact_measure}
\end{figure}

We proceed to compute the compactness measure for the hip roof house characterized by an external surface $S$ and volume $V$.  
Using equations \eqref{powierzchnia2} and \eqref{S_min}, we obtain:
\eqb
\frac{S}{S_{\min}}=\frac{V^{2/3}\left( \frac{2k + 2rk + \frac{r}{\cos(\al)}}{(rk)^{2/3}}\right)}{\frac{3\cdot (2V)^{2/3}}{(\cos(\al))^{1/3}}}=
\frac{(\cos(\al))^{1/3}\left( {2k + 2rk + \frac{r}{\cos(\al)}}\right)}{{3\cdot (2 r k)^{2/3}}}
\eqe
It should be emphasized that the ratio $\frac{S}{S_{\min}}$ computed above is indeed dimensionless and independent of the volume $V$.  
It can be factorized into $\gamma(r,k)$, as defined in equation \eqref{gamma_f}, and a term dependent on $\alpha$.  
Furthermore, $\frac{S}{S_{\min}} \geq 1$, with equality attained for the optimal parameters $r_{\min} $ and $k_{\min}$, see Fig. \ref{Compact_measure}.  
This compactness measure thus serves as a quantitative indicator of the efficiency of the hip roof house design, offering valuable insight during the design process.  
It can be employed to assess shape optimality, guide improvements, and identify directions for potential reductions in energy consumption and construction costs.  

In Sec. \ref{Case_studies}, we will apply this compactness measure to evaluate the optimality of selected existing hip roof houses in real-world settings.

\subsubsection{Second scenario: fixed footprint ratio $r$}
\label{Scenario2}

 Now, for the optimization of the hip roof house we assume a fixed footprint ratio $r$. This assumption is justified because footprint (lot coverage) is commonly constrained by planning regulations (setbacks, maximum coverage or floor-area-ratio limits), neighborhood covenants, and client requirements, so it is frequently a design parameter set prior to detailed form finding. In addition, the geometry of the building plot (shape, access, terrain and existing easements) often dictates the only feasible or economically viable footprint shapes and areas; treating  $r$ as given enables focused optimization of vertical form, orientation and envelope while holding the primary land-use parameter constant. 
As before, we also assume that $V$ and $\al$ are given.
 
To find the optimal shape of the hip roof house for fixed ratio $r$, we solve the second equation in \eqref{system} 
$$
\frac{\partial \gamma(r,k)}{\partial k} = 0
$$
with respect to $k$ treating $r$ as constant. We get that the minimum of $\gamma(r,k)$ with fixed $r$ and $\al\in(0,\pi/2)$ is attained for 
\eqb
\label{k_min_r_fixed}
k_{\min}=\frac{r}{(r+1)\cos(\al)}.
\eqe
The resulting minimal external surface of the hip roof house equals
\eqb
\label{S_min_r_fixed}
S_{\min}=V^{2/3}\gamma(r,k_{\min})=
\frac{3rV^{2/3}}{\cos(\al)\left(\frac{r^2}{\cos(\al)(r + 1)}\right)^{2/3}}.
\eqe
The other optimal parameters of the hip roof house are as follows
\eqb
\label{W_min_r_fixed}
W_{\min}=\left( \frac{V}{r\cdot k_{\min}} \right)^{1/3}=
 \left(\frac{V\cos(\al)(r + 1)}{r^2}\right)^{1/3},
\eqe
\eqb
\label{W_min_r_fixed}
L_{\min}=W_{\min}\cdot r = r
 \left(\frac{V\cos(\al)(r + 1)}{r^2}\right)^{1/3} ,
\eqe
\eqb
\label{H_min_r_fixed}
H_{\min}=W_{\min}\cdot k_{\min}=
 \frac{r}{\cos(\al)(r + 1)}
 \left(\frac{V\cos(\al)(r + 1)}{r^2}\right)^{1/3} .
\eqe

Continuing the previous example, let us take $V=400\; m^3$ and $\al=\pi/6=30^\circ$ and $r=1,5$. Then, from equations \eqref{k_min_r_fixed}--\eqref{H_min_r_fixed} we get the following optimal parameters of the house:
$k_{\min}=0,69$, $S_{\min}=274,94\;m^2 $, $W_{\min}=7,27\;m$, $L_{\min}=10,91\;m$, $H_{\min}=5,04\;m$. Note that $W_{\min}\cdot 1,5 =L_{\min}$

From a practical perspective, the parameters obtained above are directly applicable to real-world house construction. They provide design solutions that are not only feasible but also structurally efficient, cost-effective, and supportive of improved energy performance. Consequently, these results offer a solid basis for architects and engineers aiming to apply optimized hip roof house designs in practice with fixed ratio $r$.

\begin{figure}[t]
\centering
\includegraphics[width=10cm]{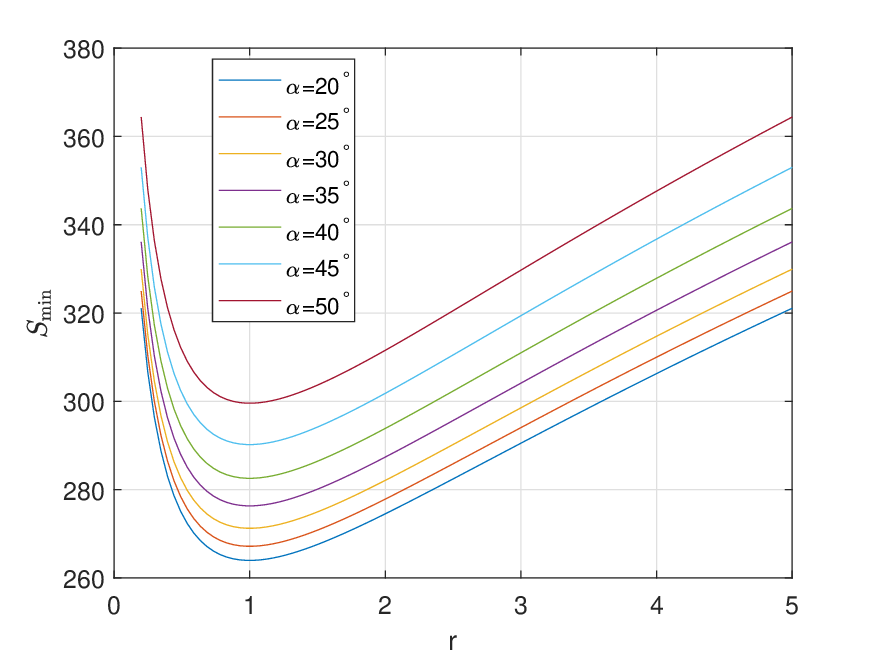}
\caption{Plots of the optimal surfaces $S_{\min}$, as functions of fixed footprint ratio $r$ for different $\al$. Here $V=400 \; m^3$.  }
\label{S_fixed_r}
\end{figure}

In Fig. \ref{S_fixed_r} we observe plots of the optimal surfaces $S_{\min}$ calculated using \eqref{S_min_r_fixed}, as functions of fixed footprint ratio $r$.  For each $\al$ we can see that the functions attain minimum for $r=1$. Here $V=400 \; m^3$.

\begin{figure}[t]
\centering
\includegraphics[width=10cm]{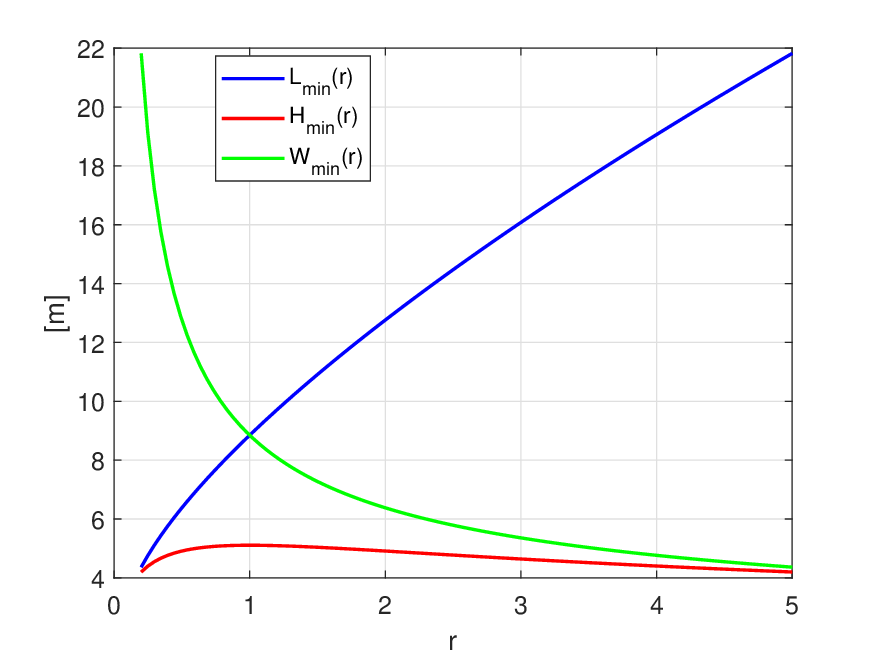}
\caption{Plots of the optimal parameters $W_{\min}$, $L_{\min}$ and $H_{\min}$ as functions of fixed footprint ratio $r$. Here $V=400 \; m^3$ and $\al=\pi/6=30^\circ$.  }
\label{Param_fixed_r}
\end{figure}

In Fig. \ref{Param_fixed_r} we present plots of the optimal parameters $W_{\min}$, $L_{\min}$ and $H_{\min}$ as functions of fixed footprint ratio $r$. Here $V=400 \; m^3$ and $\al=\pi/6=30^\circ$.
Interestingly, with increasing $r$, length $L_{\min}$ increases, $W_{\min}$ decreases, whereas $H_{\min}$ attains minimum at $r=1$.

\begin{figure}[t]
\centering
\includegraphics[width=12cm]{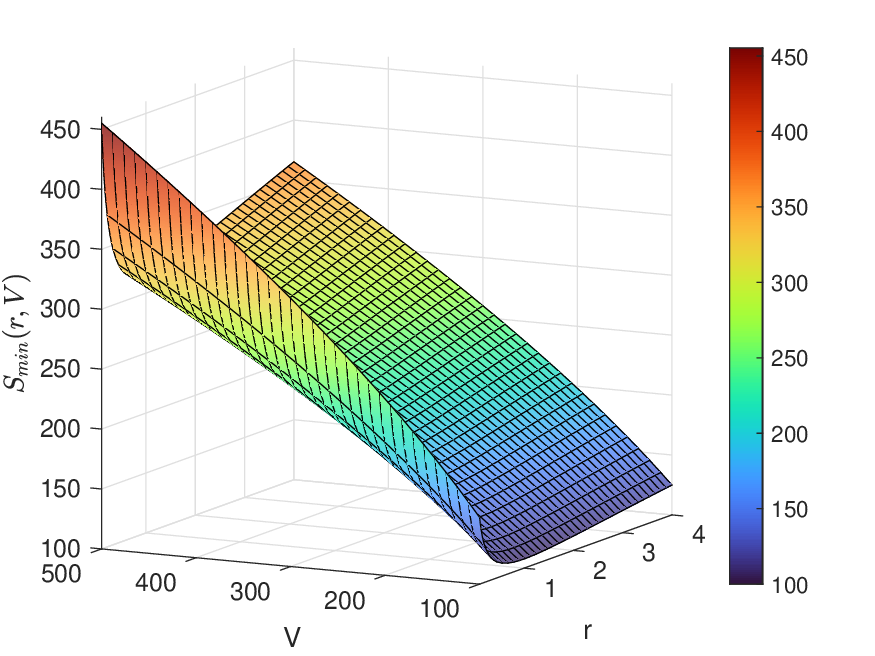}
\caption{Plot of the optimal surface $S_{\min}$, as a function of fixed footprint ratio $r$ and volume $V$. Here $\al=\pi/6=30^\circ$.}
\label{S_fixed_r_al_30}
\end{figure}

In Fig. \ref{S_fixed_r_al_30} we see a 3D plot of the optimal surface $S_{\min}$, as a function of fixed footprint ratio $r$ and volume $V$. Here $\al=\pi/6=30^\circ$. We observe that for increasing $V$ also $S_{\min}$ increases. However, for a given $V$ the surface $S_{\min}$ attains minimum at $r=1$.

\subsubsection{Third scenario: fixed slenderness ratio $k$}
\label{Scenario3}

Now, in the optimization of the hip roof house we assume a fixed slenderness ratio $k$. This assumption reflects both regulatory and functional realities of architectural and structural design. Local building codes and urban planning rules usually restrict the maximum height of residential houses in relation to their footprint, which in practice sets clear boundaries on slenderness. In many cases these limits are not negotiable, since they are intended to preserve the character of the surrounding built environment, control overshadowing and ensure adequate access to light and ventilation. At the same time, structural considerations naturally tie height to base dimension: excessively slender houses may face stability and serviceability problems, leading to higher structural costs, while overly squat proportions may be uneconomical in terms of land use and circulation efficiency. Fixing slenderness also responds to functional requirements, because internal layouts, stair and elevator cores, and the distribution of services become inefficient when the ratio departs significantly from conventional values. Moreover, energy performance and envelope-to-volume relationships are closely linked to slenderness; holding it constant allows one to isolate the impact of other geometric features such as orientation, roof form or façade articulation. In practice, the slenderness ratio is often determined early in the project by a combination of zoning constraints, structural logic and client expectations, and treating it as fixed in the optimization process mirrors this reality while ensuring that results are relevant, realistic and directly applicable to the design process.

As before, we also assume that $V$ and $\al$ are given.

To find the optimal shape of the hip roof house for fixed ratio $k$, we solve the first equation in \eqref{system} 
$$
\frac{\partial \gamma(r,k)}{\partial r} = 0
$$
with respect to $r$, treating $k$ as constant. We get that the minimum of $\gamma(r,k)$ with fixed $k$ and $\al\in(0,\pi/2)$ is attained for 
\eqb
\label{r_min_k_fixed}
r_{\min}=\frac{4k\cos(\al)}{2k\cos(\al)+1}.
\eqe
The resulting minimal external surface of the hip roof house equals
\eqb
\label{S_min_k_fixed}
S_{\min}=V^{2/3}\gamma(r_{\min},k)=
6kV^{2/3}\left(\frac{2k\cos(\al)+1}{4k^2\cos(\al)}   \right)^{2/3}.
\eqe
The other optimal parameters are as follows
\eqb
\label{W_min_k_fixed}
W_{\min}=\left( \frac{V}{r_{\min}\cdot k} \right)^{1/3}=
V^{1/3} \left(\frac{2k\cos(\al)+1}{4k^2\cos(\al)}   \right)^{1/3},
\eqe
\eqb
\label{W_min_k_fixed}
L_{\min}=W_{\min}\cdot r_{\min} = \frac{4V^{1/3}k\cos(\al)}{2k\cos(\al)+1}\left(\frac{2k\cos(\al)+1}{4k^2\cos(\al)}   \right)^{1/3},
\eqe
\eqb
\label{H_min_k_fixed}
H_{\min}=W_{\min}\cdot k= k V^{1/3} \left(\frac{2k\cos(\al)+1}{4k^2\cos(\al)}   \right)^{1/3}.
\eqe

As an example, let us take $V=400\; m^3$ and $\al=\pi/6=30^\circ$ and $k=0,5$. Then, from equations \eqref{r_min_k_fixed}--\eqref{H_min_k_fixed} we get the following optimal parameters of the house:
$r_{\min}=0,93$, $S_{\min}=271,70\;m^2 $, $W_{\min}=9,52\;m$, $L_{\min}=8,83\;m$, $H_{\min}=4,76\;m$. Note that $W_{\min}\cdot 0,5 =H_{\min}$.

From a practical standpoint, the parameters derived above can be applied to real-world house construction. They offer design solutions that are feasible, structurally efficient, cost-effective, and conducive to improved energy performance. As such, these results provide a strong foundation for architects and engineers seeking to implement optimized hip roof house designs in practice with a fixed ratio $k$.

\begin{figure}[t]
\centering
\includegraphics[width=10cm]{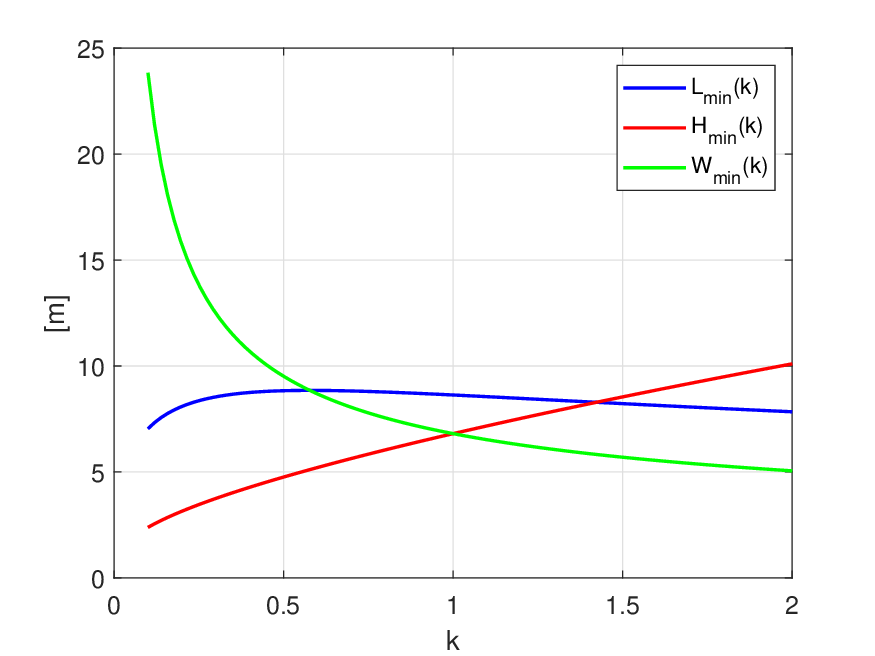}
\caption{Plots of the optimal parameters $W_{\min}$, $L_{\min}$ and $H_{\min}$ as functions of fixed slenderness ratio $k$. Here $V=400 \; m^3$ and $\al=\pi/6=30^\circ$.  }
\label{Param_fixed_k}
\end{figure}

In Fig. \ref{Param_fixed_k} can see plots of the optimal parameters $W_{\min}$, $L_{\min}$ and $H_{\min}$ as functions of fixed slenderness ratio $k$. Here $V=400 \; m^3$ and $\al=\pi/6=30^\circ$.
Interestingly, with increasing $k$, width $W_{\min}$ decreases, $H_{\min}$ increases, whereas $L_{\min}$ attains minimum at $k=\frac{1}{2\cos(\al)}=0,58$.

\begin{figure}[t]
\centering
\includegraphics[width=12cm]{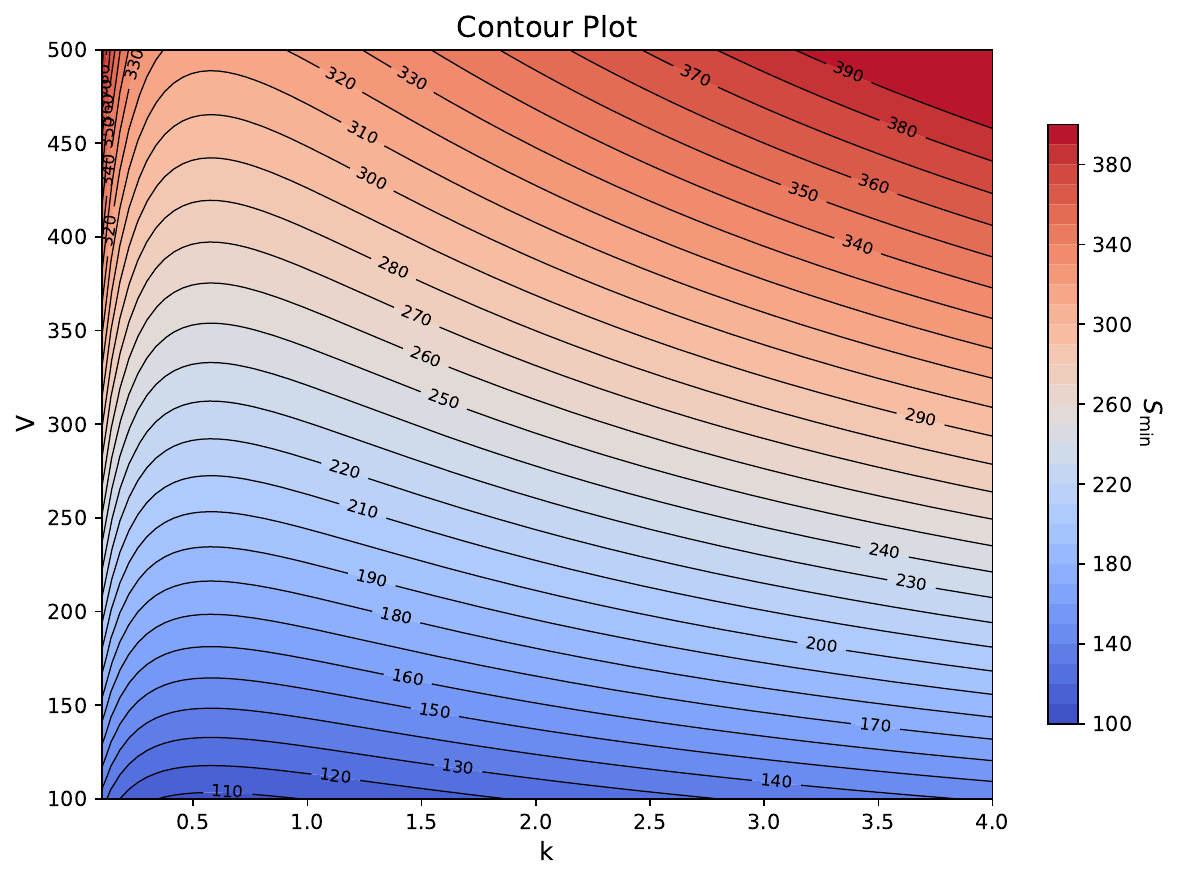}
\caption{Contour plot of the optimal surface $S_{\min}$, as a function of fixed slenderness ratio $k$ and volume $V$. Here $\al=\pi/6=30^\circ$.}
\label{contour_plot_fixed_k}
\end{figure}

In Fig. \ref{contour_plot_fixed_k} we have a contour plot of the optimal surface $S_{\min}$, as a function of fixed slenderness ratio $k$ and volume $V$. Here $\al=\pi/6=30^\circ$. We observe that for increasing $V$ also $S_{\min}$ increases. However, for a given $V$ the surface $S_{\min}$ attains minimum at $k=\frac{1}{2\cos(\al)}=0,58$. 

\subsubsection{Fourth scenario: fixed floor area}
\label{Scenario4}

Assuming that the floor area $F=W\cdot L$ is given, when optimizing the shape of a house reflects the way design constraints typically arise in practice. The required ground area is usually dictated early on by the investor’s requirements — number of rooms, circulation, and amenities—as well as by  regulations that limit lot coverage, impose setbacks, or reserve land for access and landscaping. The physical characteristics of the building plot, including its geometry, orientation, and topography, further restrict the footprint, making it more a boundary condition than a free design variable. Fixing the footprint allows the optimization to focus on aspects such as orientation, vertical arrangement, and envelope geometry, which more directly influence performance metrics like energy efficiency, daylighting, and structural behavior. It also ensures fair comparison across alternatives, since differences in results come from form rather than from changes in usable area. 

From an economic and constructability standpoint, foundations, site preparation, and utilities are all tied to a predetermined footprint, so keeping it constant avoids impractical solutions and aligns the optimization with regulatory and economic realities.

Let us now address the optimality problem for the hip roof house with given floor area. We assume that the floor area of the house $F=WL$ is fixed, and for a given 
$F$ we seek the values of $W$ and $L$ that minimize the external surface area $S$. To ensure that the problem is well-posed, the height $H$ must also be given. Otherwise, arbitrarily reducing $H$ would always decrease the surface area. Such a formulation reflects common conditions in architectural design, where the investor typically specifies not only the required floor area but also the height of the rooms. 

As before, we also assume that $V$ and $\al$ are given. In fact, $V=F\cdot H$, so knowing $F$ and $H$, we also know the value of $V$.
In such setting the external surface of the house is given by (cf. \eqref{powierzchnia}):
\begin{align}
S=S(W)=2WH+2LH+\frac{LW}{\cos(\alpha)}
= 2WH+2FH/W+\frac{F}{\cos(\alpha)}.
\label{powierzchnia3}
\end{align}
It is straightforward to find the minimum of the function $S(W)$. It attains the minimum for $W_{\min}=\sqrt{F}$. It immediately implies that also $L_{\min}=\sqrt{F}$. So, in this case, the optimal floor shape is a square, as expected.

As an example, let us take $F=100\; m^2$, $\al=\pi/6=30^\circ$ and $H=3\;m$. Then $W_{\min}=L_{\min}=\sqrt{F}=10$. The corresponding minimal surface is equal to $S_{\min}=235,47\; m^2$.

\begin{figure}[t]
\centering
\includegraphics[width=10cm]{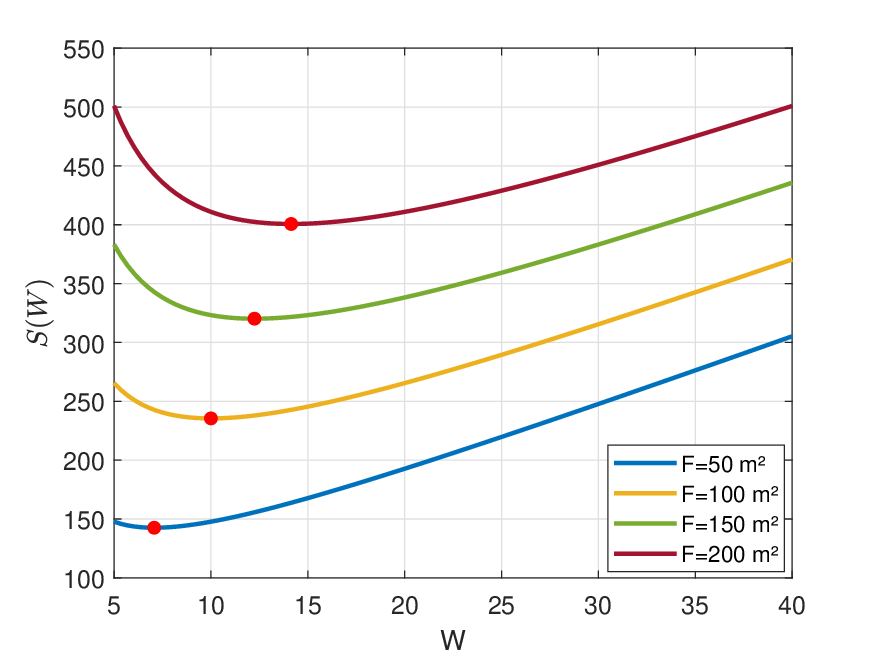}
\caption{Plots of the external surface function $S(W)$ 
for different values of the floor area $F$.
The function reaches its minimum (red dots) at $W=\sqrt{F}$. Here $\al=\pi/6=30^\circ$ and $H=3\;m^2$.}
\label{S_fixed_F}
\end{figure}

In Fig. \ref{S_fixed_F} we present plots of the external surface function $S(W)$ for different values of the floor area $F$. Indeed, the function reaches its minimum (red dots) at $W=\sqrt{F}$.

\subsubsection{Fifth scenario: height from a fixed range}
\label{Scenario5}

Assuming that the height of the house is constrained within a predefined range is a realistic and necessary assumption when approaching the problem of shape optimization. In almost all practical contexts, in the majority of developed countries, the maximum and minimum allowable heights of residential buildings are strongly influenced by zoning regulations, neighborhood planning rules, and building codes that aim to preserve the character of the built environment, protect access to sunlight for neighboring plots, and ensure fire safety or structural stability. 

At the same time, the functional program of a house imposes a lower bound on height, since a certain number of room heights are required to meet comfort and usability standards. 

Structural and economic considerations further reinforce the need for height limits: excessively tall single-family houses within a small footprint quickly become uneconomical due to the costs of vertical circulation, load-bearing systems, and additional safety requirements, while overly low proportions may lead to inefficient use of space, cramped interiors, and difficulty integrating necessary technical systems. 

Energy performance also depends on height, as very tall or very shallow building volumes alter envelope-to-volume ratios and consequently heating and cooling demands. By defining a feasible range for house height, the optimization process remains grounded in realistic design parameters, ensuring that the solutions are both technically viable and compliant with external constraints. This approach reduces the risk of producing theoretically optimal but practically unbuildable shapes, and it allows the analysis to concentrate on meaningful geometric and spatial requirements that reflect regulatory, functional, and economic realities.

So, in this scenario, we assume that the height $H$ belongs to a predefined interval, say $H\in[a,b]$. As before, we also assume that $V$ and $\al$ are given.
In such setting we will find optimal shape of the house, which will minimize the external surface $S$ given by \eqref{powierzchnia}. 
For this purpose we will use the method of Lagrange multipliers 
with Karush–Kuhn–Tucker (KKT) conditions \cite{boyd2004convex}.

We want to find the optimal parameters $(W^*, H^*)$, which will minimize the external surface (recall formula \eqref{powierzchnia})
\[
S=S(W,H) \;=\; 2WH + \frac{2V}{W} + \frac{V}{H\cos(\alpha)}
\]
subject to constrains
\[
a \leq H \leq b.
\]
Here, $V, a, b$, and $\alpha$ are positive constants, with the assumptions that $0< a < b$ and $\cos(\alpha) > 0$. We are looking for solutions where $W > 0$ and $H > 0$.

First, we rewrite the constraints in the standard inequality form using auxiliary functions $g_i(W,H),\;i=1,2$:
\begin{align*}
H \geq a &\implies g_1(W,H) = a - H \le 0 \\
H \leq b &\implies g_2(W,H) = H - b \le 0
\end{align*}

The Lagrangian for this optimization problem $L(W, H, \mu_1, \mu_2)$ is defined as:
\[
L(W, H, \mu_1, \mu_2) = S(W,H) + \mu_1 g_1(W,H) + \mu_2 g_2(W,H),
\]
where $\mu_1$ and $\mu_2$ are the Lagrange multipliers.

Substituting  functions $S$, $g_1$ and $g_2$, we get:
\[
L(W, H, \mu_1, \mu_2) = 2WH + \frac{2V}{W} + \frac{V}{H\cos(\alpha)} + \mu_1(a - H) + \mu_2(H - b).
\]

For a point $(W^*, H^*)$ to be an optimal solution, it must satisfy the following four KKT conditions \cite{boyd2004convex}:

\begin{enumerate}
    \item[1.] \textbf{Stationarity:} The gradient of the Lagrangian with respect to the variables is zero
    \begin{align}
        \frac{\partial L}{\partial W} &= 2H - \frac{2V}{W^2} = 0 \label{eq:stat1} \\
        \frac{\partial L}{\partial H} &= 2W - \frac{V}{H^2\cos(\alpha)} - \mu_1 + \mu_2 = 0 \label{eq:stat2}
    \end{align}

    \item[2.] \textbf{Primal Feasibility:} The constraints must be satisfied
    \begin{align}
        a - H &\le 0 \implies H \ge a \label{eq:primal1} \\
        H - b &\le 0 \implies H \le b \label{eq:primal2}
    \end{align}

    \item[3.] \textbf{Dual Feasibility:} The Lagrange multipliers must be non-negative
    \begin{align}
        \mu_1 &\ge 0 \label{eq:dual1} \\
        \mu_2 &\ge 0 \label{eq:dual2}
    \end{align}

    \item[4.] \textbf{Complementary Slackness:}
    \begin{align}
        \mu_1(a - H) &= 0 \label{eq:slack1} \\
        \mu_2(H - b) &= 0 \label{eq:slack2}
    \end{align}
\end{enumerate}

From the first stationarity condition \eqref{eq:stat1}, we have:
\begin{equation}
2H = \frac{2V}{W^2} \implies H = \frac{V}{W^2} \implies W = \sqrt{\frac{V}{H}}. \label{eq:W_H_relation}
\end{equation}

We will find the solution by considering all four cases for the complementary slackness condition:
\paragraph*{Case 1: No constraints are active ($\mu_1 = 0$ and $\mu_2 = 0$)}

This case corresponds to the solution being in the interior of the feasible region, i.e., $a \leq H \leq b$. The second stationarity condition \eqref{eq:stat2} becomes:
\begin{equation}
2W - \frac{V}{H^2\cos(\alpha)} = 0 \implies 2W = \frac{V}{H^2\cos(\alpha)} \label{eq:case1_stat2}
\end{equation}
Substitute \eqref{eq:W_H_relation} into \eqref{eq:case1_stat2}:
\begin{align*}
2W &= \frac{V}{\left(\frac{V}{W^2}\right)^2\cos(\alpha)} = \frac{V W^4}{V^2\cos(\alpha)} = \frac{W^4}{V\cos(\alpha)} \\
2 &= \frac{W^3}{V\cos(\alpha)} \quad (\text{since } W>0) \\
W^3 &= 2V\cos(\alpha) \implies W^* = (2V\cos(\alpha))^{1/3}
\end{align*}
Now, we find the corresponding $H^*$:
\begin{align*}
H^* = \frac{V}{(W^*)^2} = \frac{V}{(2V\cos(\alpha))^{2/3}} = \frac{V}{V^{2/3}(2\cos(\alpha))^{2/3}} = V^{1/3} \left(\frac{1}{4\cos^2(\alpha)}\right)^{1/3} = \left(\frac{V}{4\cos^2(\alpha)}\right)^{1/3}
\end{align*}
This solution is valid if and only if $a < H^* < b$.

\paragraph*{Case 2: The first constraint is active ($\mu_1 > 0$ and $\mu_2 = 0$)}

From the slackness condition \eqref{eq:slack1}, $a - H = 0 \implies H = a$.
Using \eqref{eq:W_H_relation}, we find $W$:
\[
W = \sqrt{\frac{V}{a}}
\]
From the second stationarity condition \eqref{eq:stat2}, we solve for $\mu_1$:
\[
\mu_1 = 2W - \frac{V}{a^2\cos(\alpha)} = 2\sqrt{\frac{V}{a}} - \frac{V}{a^2\cos(\alpha)}
\]
The dual feasibility condition $\mu_1 \ge 0$ implies:
\begin{align*}
2\sqrt{\frac{V}{a}} - \frac{V}{a^2\cos(\alpha)} \ge 0 &\implies 2\sqrt{\frac{V}{a}} \ge \frac{V}{a^2\cos(\alpha)} \\
&\implies 2 \ge \frac{\sqrt{V}}{a^{3/2}\cos(\alpha)} \\
&\implies 2a^{3/2}\cos(\alpha) \ge \sqrt{V} \\
&\implies 4a^3\cos^2(\alpha) \ge V \\
&\implies a^3 \ge \frac{V}{4\cos^2(\alpha)} \\
&\implies a \ge \left(\frac{V}{4\cos^2(\alpha)}\right)^{1/3}
\end{align*}
This case occurs when the unconstrained optimum $H^*$ is less than or equal to $a$.

\paragraph*{Case 3: The second constraint is active ($\mu_1 = 0$ and $\mu_2 > 0$)}

From the slackness condition \eqref{eq:slack2}, $H - b = 0 \implies H = b$.
Using \eqref{eq:W_H_relation}, we find $W$:
\[
W = \sqrt{\frac{V}{b}}
\]
From the second stationarity condition \eqref{eq:stat2}, we solve for $\mu_2$:
\[
\mu_2 = \frac{V}{b^2\cos(\alpha)} - 2W = \frac{V}{b^2\cos(\alpha)} - 2\sqrt{\frac{V}{b}}
\]
The dual feasibility condition $\mu_2 \ge 0$ implies:
\begin{align*}
\frac{V}{b^2\cos(\alpha)} - 2\sqrt{\frac{V}{b}} \ge 0 &\implies \frac{V}{b^2\cos(\alpha)} \ge 2\sqrt{\frac{V}{b}} \\
&\implies \frac{\sqrt{V}}{b^{3/2}\cos(\alpha)} \ge 2 \\
&\implies \sqrt{V} \ge 2b^{3/2}\cos(\alpha) \\
&\implies V \ge 4b^3\cos^2(\alpha) \\
&\implies \frac{V}{4\cos^2(\alpha)} \ge b^3 \\
&\implies \left(\frac{V}{4\cos^2(\alpha)}\right)^{1/3} \ge b
\end{align*}
This case occurs when the unconstrained optimum $H^*$ is greater than or equal to $b$.

\paragraph*{Case 4: Both constraints are active ($\mu_1 > 0$ and $\mu_2 > 0$)}

This would imply $H = a$ and $H = b$. Since we are given that $a < b$, this case is impossible.

\paragraph*{Summary of the Solution}

Let's define the critical value $H_{crit} = \left(\frac{V}{4\cos^2(\alpha)}\right)^{1/3}$. The optimal solution $(W^*, H^*)$ depends on the position of $H_{crit}$ relative to the interval $[a, b]$.
This can be written compactly in the following way:
\eqb
(W^*, H^*) =
\begin{cases}
    \left(\sqrt{\frac{V}{a}}, a\right) & \text{if } \left(\frac{V}{4\cos^2(\alpha)}\right)^{1/3} \le a \\
    \\
    \left((2V\cos(\alpha))^{1/3}, \left(\frac{V}{4\cos^2(\alpha)}\right)^{1/3}\right) & \text{if } a < \left(\frac{V}{4\cos^2(\alpha)}\right)^{1/3} < b \\
    \\
    \left(\sqrt{\frac{V}{b}}, b\right) & \text{if } \left(\frac{V}{4\cos^2(\alpha)}\right)^{1/3} \ge b
\end{cases}
\label{solution_KKT}
\eqe
Consequently, the optimal length $L^*=\frac{V}{W^* H^*}$.
Comparing the above results with eqs. \eqref{W_min}-\eqref{H_min} we can see that when the critical value $H_{crit}$ belongs to the interval $(a,b)$, the optimal solution $(W^*, H^*)$ is equal to the one derived in \eqref{W_min}-\eqref{H_min}. Otherwise, the solutions differ. 

\begin{figure}[t]
\centering
\includegraphics[width=14cm]{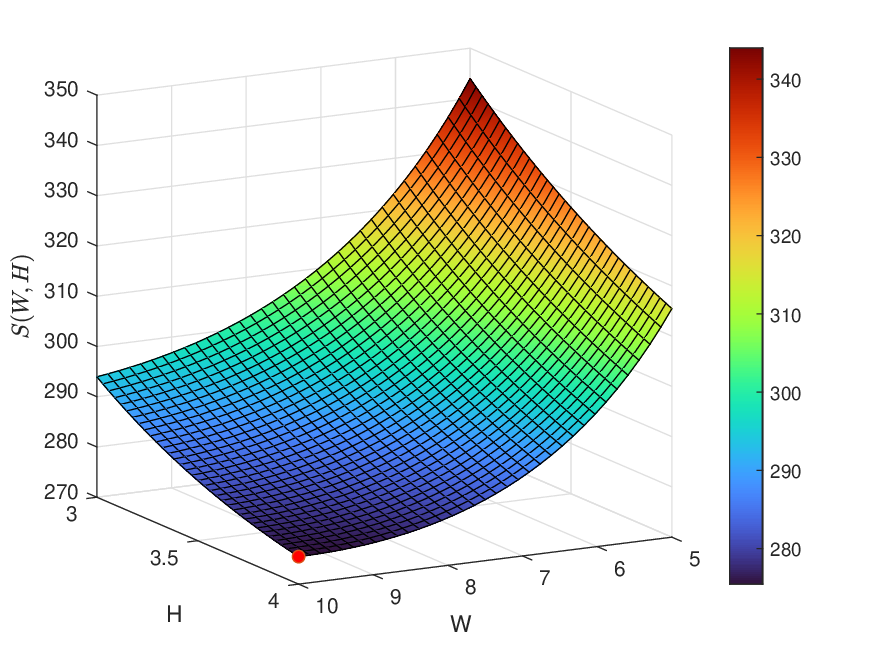}
\caption{Graph of the external surface $S(W, H)$ with parameters $V = 400\;m^3$, $\al=\pi/6=30^\circ$, $a=3$, $b=4$ . The red dot indicates the location of the
minimal surface derived using formula \eqref{solution_KKT}
for the case $H_{crit}\geq b$.}
\label{S_H_a_b}
\end{figure}

Let us consider the following example. Take $V=400\; m^3$, $\al=\pi/6=30^\circ$, $a=3\; m$ and $b=4\;m$. Then, the critical value $H_{crit}=5,11$. So, $H_{crit}\geq b$.
In this case, we get from \eqref{solution_KKT} that $(W^*, H^*)=\left(\sqrt{\frac{V}{b}}, b\right)=
\left(10\;,\;4\right)$, and the minimal surface is $S_{\min}=275,47\;m^2$, see Fig. \ref{S_H_a_b}. 

\begin{figure}[t]
\centering
\includegraphics[width=14cm]{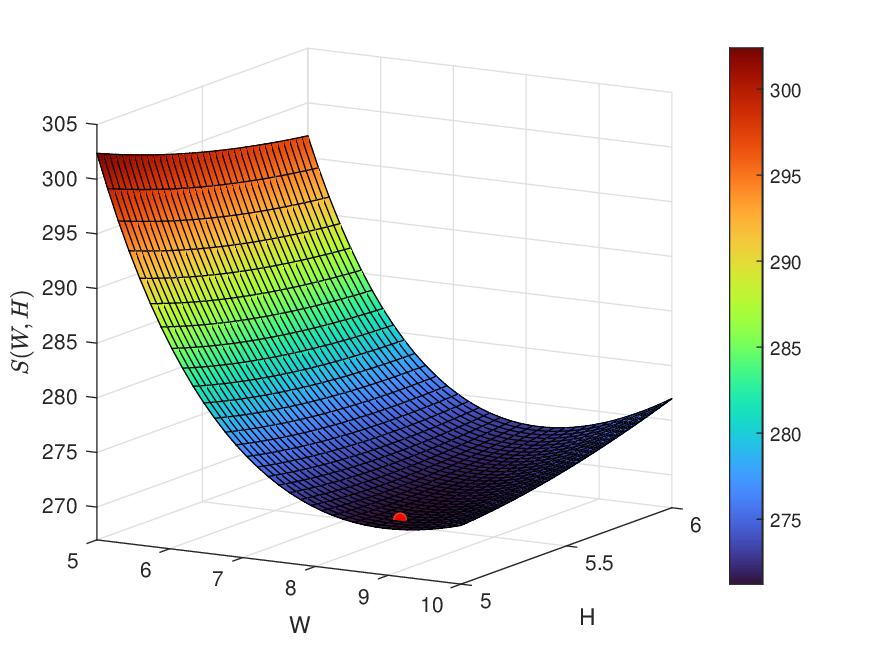}
\caption{Plots of the external surface $S(W, H)$ with parameters $V = 400\;m^3$, $\al=\pi/6=30^\circ$, $a=5$, $b=6$ . The red dot indicates the location of the
minimal surface derived using formula \eqref{solution_KKT}
for the case $H_{crit}\in (a,b)$.}
\label{S_H_a_b2}
\end{figure}

However, if we take $a=5\; m$ and $b=6\;m$ and the other parameters as above, then $H_{crit}\in (a,b)$.
In this case, we get from \eqref{solution_KKT} that $(W^*, H^*)= \left((2V\cos(\alpha))^{1/3}, \left(\frac{V}{4\cos^2(\alpha)}\right)^{1/3}\right)=(8,85\; ; \;5,11)$, and the minimal surface is $S_{\min}=271,23\;m^2$, see Fig. \ref{S_H_a_b2}.

\section{Case studies}
\label{Case_studies}
 
\begin{figure}[t]
\centering
\includegraphics[width=14cm]{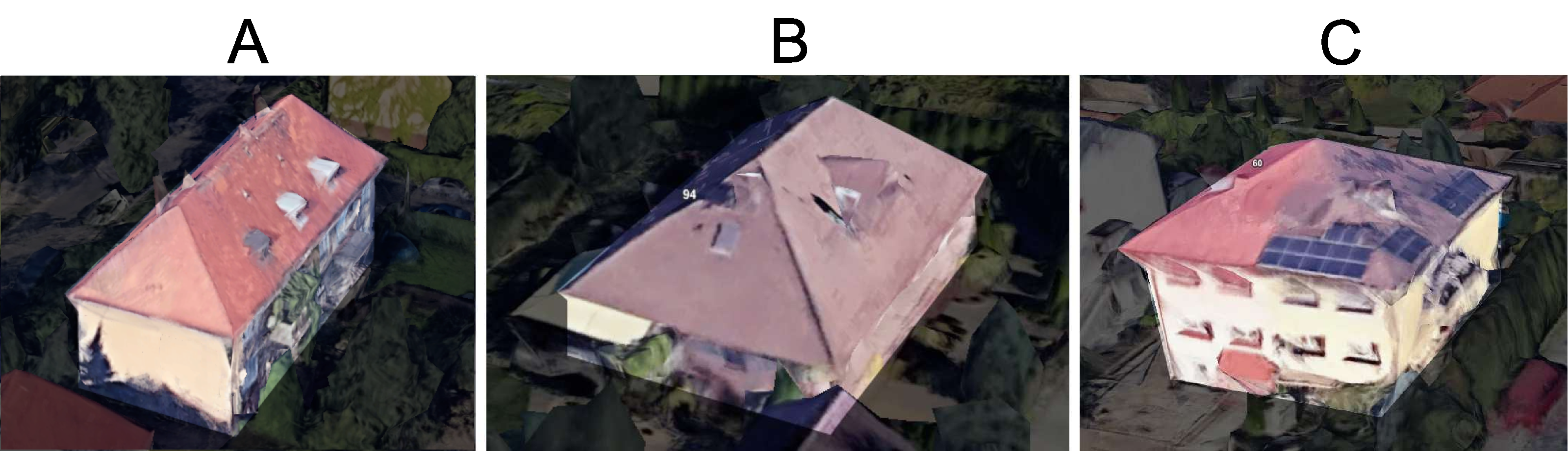}
\caption{Images of three hip roof houses analyzed as case studies using the results obtained
in previous sections. Images source: Google Earth. }
\label{cases}
\end{figure}

To highlight the practical relevance of the theoretical optimality results established in the preceding section, we present an application in the form of three case studies. Specifically, we investigate three recently constructed hip roof houses, labeled A, B, and C (see Fig. \ref{cases}). These buildings, located in Wrocław, Poland, were chosen deliberately to capture a representative variety of geometric forms, proportions, and overall scales. The aim is to demonstrate how the derived formulas can be employed to analyze real structures of different dimensions and configurations. For each house, both the visual appearance and the essential geometric parameters were obtained using publicly available satellite image and measurement tools provided by Google Earth. This approach ensures that the case studies are not only illustrative but also reproducible, since the case studies can be extended to other architectural examples.

\paragraph{Analysis of House A}\ \\
The parameters of House A are the following: \\
width $W=10,9\;m$,\\
length $L=26,7\;m$,\\
height $H=7,2\;m$,\\
roof slope angle $\al=50^\circ$.\\
This gives the volume $V = 2095,4\; m^3$ and 
external surface $S = 994,2 m^2$.

Let us first apply the optimality method from Scenario (\ref{Scenario1}).
For a given fixed $V = 2095,4\; m^3$ and $\al=50^\circ$, from Eqs. \eqref{S_min}--\eqref{H_min} we get the following optimal parameters of the house:
$S_{\min} = 903,58 \; m^2$, 
$W_{\min} = 13,9 \;m$,
$L_{\min} = 13,9 \; m$,
$H_{\min} = 10,8 \; m$.
Consequently, the compactness measure equals $\frac{S}{S_{\min}}= 1,10$.
The above optimization procedure would reduce the external surface by $S-S_{\min}=90,7 \; m^2$. It would result in significant savings in construction costs and energy consumption of the house.

Let us now apply the optimality method from Scenario (\ref{Scenario2}).
For a given fixed $r=L/W=2,4$, from Eqs. \eqref{S_min_r_fixed}--\eqref{H_min_r_fixed} we get the following optimal parameters of the house:
$S_{\min} = 964,0 \; m^2$,
$W_{\min} = 9,2 \;m$,
$L_{\min} = 22,5 \; m$,
$H_{\min} = 10,1 \; m$.
Consequently, the compactness measure equals $\frac{S}{S_{\min}}= 1,03$.
This optimization lowers the external surface by 
$S-S_{\min}=30,2 \; m^2$, resulting in considerable cost savings and improved energy efficiency of the house.

Let us apply the optimality method from Scenario (\ref{Scenario3}).
For a given fixed $k=H/W=0,66$, from Eqs. \eqref{S_min_k_fixed}--\eqref{H_min_k_fixed} we get the following optimal parameters of the house:
$S_{\min} = 905,6 \; m^2$,
$W_{\min} = 15,1 \;m$,
$L_{\min} = 13,9 \; m$,
$H_{\min} = 10,0 \; m$.
Consequently, the compactness measure equals $\frac{S}{S_{\min}}= 1,10$.
This optimization can lower the external surface by 
$S-S_{\min}=88,61 \; m^2$, resulting in impressive cost savings and improved energy efficiency of the house.

Now, we apply the optimality method from Scenario (\ref{Scenario4}).
For a given fixed floor area $F=W\cdot L=291,0\; m^2$, we get the following optimal parameters of the house:
$S_{\min} = 944,1 \; m^2$,
$W_{\min} = 17,1 \;m$,
$L_{\min} = 17,1 \; m$,
Consequently, the compactness measure equals $\frac{S}{S_{\min}}= 1,05$.
This optimization lowers the external surface by 
$S-S_{\min}=50,1 \; m^2$, resulting in considerable cost savings and improved energy efficiency of the house.

\paragraph{Analysis of House B}\ \\
The parameters of House B are as follows: \\
width $W=9,5\;m$,\\
length $L=16,7\;m$,\\
height $H=2,6\;m$,\\
roof slope angle $\al=30^\circ$.\\
This gives the volume $V = 412,5\; m^3$ and 
external surface $S = 319,4 m^2$.

Let us apply the optimality method from Scenario (\ref{Scenario1}).
For a given fixed $V = 412,5\; m^3$ and $\al=30^\circ$, from Eqs. \eqref{S_min}--\eqref{H_min} we get the following optimal parameters of the house:
$S_{\min} = 276,8 \; m^2$,
$W_{\min} = 8,9 \;m$,
$L_{\min} = 8,9 \; m$,
$H_{\min} = 7,7 \; m$.
Consequently, the compactness measure equals $\frac{S}{S_{\min}}= 1,15$.
The proposed optimization reduces the external surface area by 
$S-S_{\min}=42,7 \; m^2$, leading to substantial savings in both construction costs and the house’s energy consumption.

Now, we apply the optimality method from Scenario (\ref{Scenario2}).
For a given fixed $r=L/W=1,7$, from Eqs. \eqref{S_min_r_fixed}--\eqref{H_min_r_fixed} we get the following optimal parameters:
$S_{\min} = 284,2 \; m^2$,
$W_{\min} = 6,8 \;m$,
$L_{\min} = 12,0 \; m$,
$H_{\min} = 5,0 \; m$.
Consequently, the compactness measure equals $\frac{S}{S_{\min}}= 1,12$.
Applying this optimization decreases the external surface area by $S-S_{\min}=35,2 \; m^2$, which translates into notable reductions in construction expenses and energy use.

Now, let us apply the optimality method from Scenario (\ref{Scenario3}).
For a given fixed $k=H/W=0,27$, from Eqs. \eqref{S_min_k_fixed}--\eqref{H_min_k_fixed} we obtain the following optimal parameters:
$S_{\min} = 289,7 \; m^2$,
$W_{\min} = 13,3 \;m$,
$L_{\min} = 8,5 \; m$,
$H_{\min} = 3,6 \; m$.
Consequently, the compactness measure equals $\frac{S}{S_{\min}}= 1,10$.
This optimization decreases the external surface area by $S-S_{\min}=29,7 \; m^2$, which translates into notable reductions in construction expenses and energy use.

Now, we apply the optimality method from Scenario (\ref{Scenario4}).
For a given fixed floor area $F=W\cdot L=158,6\; m^2$, we get the following optimal parameters of the house:
$S_{\min} = 314,2 \; m^2$,
$W_{\min} = 12,6 \;m$,
$L_{\min} = 12,6 \; m$,
Consequently, the compactness measure equals $\frac{S}{S_{\min}}= 1,02$.
This optimization lowers the external surface by 
$S-S_{\min}=5,2 \; m^2$, resulting in small cost savings.

\paragraph{Analysis of House C}\ \\
The parameters of House C are the following: \\
width $W=12,5\;m$,\\
length $L=12,5\;m$,\\
height $H=7,9\;m$,\\
roof slope angle $\al=35^\circ$.\\
This gives the volume $V = 1234,4\; m^3$ and 
external surface $S = 585,7,4 m^2$.

Now, we apply the optimality method from Scenario (\ref{Scenario1}).
For a given fixed $V = 1234,4\; m^3$ and $\al=35^\circ$, from Eqs. \eqref{S_min}--\eqref{H_min} we get the following optimal parameters:
$S_{\min} = 585,7 \; m^2$,
$W_{\min} = 12,6 \;m$,
$L_{\min} = 12,6 \; m$,
$H_{\min} = 7,9 \; m$.
Consequently, the compactness measure equals $\frac{S}{S_{\min}}= 1,00$.
The proposed optimization reduces the external surface area by 
$S-S_{\min}=0,08 \; m^2$. This result shows that House C is very close to optimal.

Now, we apply the optimality method from Scenario (\ref{Scenario2}).
For a given fixed $r=L/W=1$, from Eqs. \eqref{S_min_r_fixed}--\eqref{H_min_r_fixed} we get the following optimal parameters:
$S_{\min} = 585,7 \; m^2$,
$W_{\min} = 12,6 \;m$,
$L_{\min} = 12,6 \; m$,
$H_{\min} = 7,7 \; m$.
Consequently, the compactness measure equals $\frac{S}{S_{\min}}= 1,00$.
This optimization reduces the external surface area by 
$S-S_{\min}=0,08 \; m^2$. This result shows again that House C is very close to optimal.

Now, we apply the optimality method from Scenario (\ref{Scenario3}).
For a given fixed $k=H/W=0,63$, from Eqs. \eqref{S_min_k_fixed}--\eqref{H_min_k_fixed} we obtain the following optimal parameters:
$S_{\min} = 585,7 \; m^2$,
$W_{\min} = 12,4 \;m$,
$L_{\min} = 12,6 \; m$,
$H_{\min} = 7,9 \; m$.
Consequently, the compactness measure equals $\frac{S}{S_{\min}}= 1,00$.
This optimization procedure reduces the external surface area by 
$S-S_{\min}=0,02 \; m^2$. This result shows again that House C is very close to optimal.

Now, we apply the optimality method from Scenario (\ref{Scenario4}).
For a given fixed floor area $F=W\cdot L=156,2\; m^2$, we get the following optimal parameters of the house:
$S_{\min} = 585,7 \; m^2$,
$W_{\min} = 12,5 \;m$,
$L_{\min} = 12,5 \; m$,
Consequently, the compactness measure equals $\frac{S}{S_{\min}}= 1,00$.
The external surface does not change, showing again that House C is optimal.

Summary of the optimization procedure from Scenario (\ref{Scenario1}) (fixed volume), applied to the three case studied hip roof houses, can be found in Tab. \ref{tab:volume}. 
For the summary of the optimization procedure from Scenario (\ref{Scenario2}) with fixed ratio $r$, see Tab. \ref{tab:fixed_r}. 
Summary of the optimization procedure from Scenario (\ref{Scenario3}) for fixed ratio $k$ is presented in Tab. \ref{tab:fixed_k}. 
Brief summary of the optimization procedure from Scenario (\ref{Scenario4}) for the case of fixed floor area is presented in Tab. \ref{tab:area}.

In the above case studies, we did not use the results from Scenario (\ref{Scenario5}), since it requires a predefined range for the height $H$. This optimization procedure can be applied at the design stage, in accordance with the investor's guidelines.

\begin{table}[h!]
\scriptsize
\centering
\caption{Summary of the optimization procedure from Scenario (\ref{Scenario1}) (fixed volume) for the three case studied hip roof houses.}
\label{tab:volume}
\begin{tabular}{|l|c|c|c|c|c|c|c|c|c|c|c|c|}
\hline
& \multicolumn{6}{c|}{\textbf{Real parameters}} & \multicolumn{4}{c|}{\textbf{Optimal parameters}} & \multicolumn{2}{c|}{\textbf{Comparison}} \\
\hline
& \makecell{$W$\\$[m]$} & \makecell{$L$\\$[m]$} & \makecell{$H$\\$[m]$ } & \makecell{$\al$\\$[^\circ]$ } & \makecell{$V$\\ $[m^3]$} & \makecell{$S$\\ $[m^2]$} & \makecell{$W_{\min}$\\ $[m]$} & \makecell{$L_{\min}$\\$[m]$} & \makecell{$H_{\min}$\\ $[m]$} & \makecell{$S_{\min}$\\ $[m^2]$} & \makecell{$S/S_{\min}$} & \makecell{$S-S_{\min}$\\ $[m^2]$} \\
\hline
\textbf{House A} &10,9&26,7&7,2&50&2095,4&994,2&13,9 &13,9 &10,8  &903,5  &1,10  &90,7   \\
\hline
\textbf{House B} &9,5 &16,7  &2,6  &30  &412,5  &319,4  & 8,9 &8,9  &5,1  &276,8  &1,15  &42,6    \\
\hline
\textbf{House C} & 12,5 &12,5  &7,9  &35  &1234,4  &585,7  & 12,6 &12,6  &7,7  &585,7  &1,00  &0,08    \\
\hline
\end{tabular}
\end{table}

\begin{table}[h!]
\scriptsize
\centering
\caption{Summary of the optimization procedure from Scenario (\ref{Scenario2}) (fixed ratio $r$) for the three case studied hip roof houses.}
\label{tab:fixed_r}
\begin{tabular}{|l|c|c|c|c|c|c|c|c|c|c|c|c|}
\hline
& \multicolumn{6}{c|}{\textbf{Real parameters}} & \multicolumn{4}{c|}{\textbf{Optimal parameters}} & \multicolumn{2}{c|}{\textbf{Comparison}} \\
\hline
& \makecell{$W$\\$[m]$} & \makecell{$L$\\$[m]$} & \makecell{$H$\\$[m]$ } & \makecell{$\al$\\$[^\circ]$ } & \makecell{$r$\\} & \makecell{$S$\\ $[m^2]$} & \makecell{$W_{\min}$\\ $[m]$} & \makecell{$L_{\min}$\\$[m]$} & \makecell{$H_{\min}$\\ $[m]$} & \makecell{$S_{\min}$\\ $[m^2]$} & \makecell{$S/S_{\min}$} & \makecell{$S-S_{\min}$\\ $[m^2]$} \\
\hline
\textbf{House A} &10,9&26,7&7,2&50& 2,4  &994,2& 9,2  &22,5   &10,1   & 964,0  &1,03   &30,2      \\
\hline
\textbf{House B} &9,5 &16,7  &2,6  &30  &  1,7  &319,4  &  6,8  & 12,0  & 5,0  & 284,2  &  1,12 & 35,2       \\
\hline
\textbf{House C} & 12,5 &12,5  &7,9  &35  &  1 &585,7  &  12,6 & 12,6  & 7,7  & 585,7  & 1,00  &  0,08       \\
\hline
\end{tabular}
\end{table}

\begin{table}[h!]
\scriptsize
\centering
\caption{Summary of the optimization procedure from Scenario (\ref{Scenario3}) (fixed ratio $k$) for the three case studied hip roof houses.}
\label{tab:fixed_k}
\begin{tabular}{|l|c|c|c|c|c|c|c|c|c|c|c|c|}
\hline
& \multicolumn{6}{c|}{\textbf{Real parameters}} & \multicolumn{4}{c|}{\textbf{Optimal parameters}} & \multicolumn{2}{c|}{\textbf{Comparison}} \\
\hline
& \makecell{$W$\\$[m]$} & \makecell{$L$\\$[m]$} & \makecell{$H$\\$[m]$ } & \makecell{$\al$\\$[^\circ]$ } & \makecell{$k$\\} & \makecell{$S$\\ $[m^2]$} & \makecell{$W_{\min}$\\ $[m]$} & \makecell{$L_{\min}$\\$[m]$} & \makecell{$H_{\min}$\\ $[m]$} & \makecell{$S_{\min}$\\ $[m^2]$} & \makecell{$S/S_{\min}$} & \makecell{$S-S_{\min}$\\ $[m^2]$} \\
\hline
\textbf{House A} &10,9&26,7&7,2&50& 0,66  &994,2  &  15,1 & 13,9  & 10,0  & 905,6   &1,10   & 88,61     \\
\hline
\textbf{House B} &9,5 &16,7  &2,6  &30  &  0,27  &319,4  &13,3   &8,5   & 3,6  &   289,7 & 1,10  &   29,7    \\
\hline
\textbf{House C} & 12,5 &12,5  &7,9  &35  & 0,63  &585,7  & 12,4 & 12,6  & 7,9  &  585,7  & 1,00  &  0,02    \\
\hline
\end{tabular}
\end{table}

\begin{table}[h!]
\scriptsize
\centering
\caption{Summary of the optimization procedure from Scenario (\ref{Scenario4}) (fixed floor area) for the three case studied houses.}
\label{tab:area}
\begin{tabular}{|l|c|c|c|c|c|c|c|c|c|c|c|}
\hline
& \multicolumn{6}{c|}{\textbf{Real parameters}} & \multicolumn{3}{c|}{\textbf{Optimal parameters}} & \multicolumn{2}{c|}{\textbf{Comparison}} \\
\hline
& \makecell{$W$\\$[m]$} & \makecell{$L$\\$[m]$} & \makecell{$H$\\$[m]$ } & \makecell{$\al$\\$[^\circ]$ } & \makecell{$F$\\ $[m^2]$} & \makecell{$S$\\ $[m^2]$} & \makecell{$W_{\min}$\\ $[m]$} & \makecell{$L_{\min}$\\$[m]$} & \makecell{$S_{\min}$\\ $[m^2]$} & \makecell{$S/S_{\min}$} & \makecell{$S-S_{\min}$\\ $[m^2]$} \\
\hline
\textbf{House A} &10,9&26,7&7,2&50&  291,0 &994,2  & 17,1   & 17,1   &  944,1  & 1,05   & 50,1  \\
\hline
\textbf{House B} &9,5 &16,7  &2,6  &30  &  158,6  &319,4& 12,6   &12,6    & 314,2   & 1,02   &  5,2    \\
\hline
\textbf{House C} & 12,5 &12,5  &7,9  &35  & 156,2  &585,7 & 12,5   & 12,5   & 585,7   &  1,00  & 0  \\
\hline
\end{tabular}
\end{table}

\paragraph{Discussion on case studies} \ \\
The comparative analysis of Houses A, B, and C reveals clear differences in how geometric proportions affect the efficiency of hip roof houses. Although all three buildings share the same roof typology, their deviations from the theoretical optimum vary considerably, demonstrating the importance of proportion in architectural design.

House A is defined by an elongated rectangular footprint with a footprint ratio $r = 2,4$ and a steep roof slope of $\alpha = 50^\circ$. While its large volume of $2095,4 \,\mathrm{m}^3$ provides significant interior space, the corresponding external surface of $994,2 \,\mathrm{m}^2$ is far from optimal. The compactness measure, which ranges from $1,03$ in the case of a fixed footprint ratio to $1,10$ when evaluated under fixed volume or slenderness conditions, indicates an excess envelope of up to $90 \,\mathrm{m}^2$. This surplus translates into higher construction costs and increased thermal losses. The analysis suggests that the elongated plan is the primary factor behind this inefficiency, as designs constrained to a rectangular footprint tend to move away from the compactness inherent in square-based forms.

House B presents a different situation. With $W = 9,5 \,\mathrm{m}$, $L = 16,7 \,\mathrm{m}$, and a very low height of $H = 2,6 \,\mathrm{m}$, the building has a slenderness ratio of only $k = 0,27$. The resulting surface area of $319,4 \,\mathrm{m}^2$ exceeds the theoretical minimum up to $43 \,\mathrm{m}^2$, depending on the optimization scenario. The compactness measure is correspondingly high, reaching $1,15$ in the fixed-volume case. The analysis reveals that the low vertical dimension is the principal source of inefficiency: by compressing the height, the building increases its surface-to-volume ratio, thereby undermining both thermal performance and material efficiency. Even moderate increases in height would bring the design significantly closer to the optimum.

House C, by contrast, is an almost exact match to the theoretical predictions. With a square footprint of $12,5 \,\mathrm{m} \times 12,5 \,\mathrm{m}$, a height of $7,9 \,\mathrm{m}$, and a roof slope angle of $35^\circ$, the resulting ratios $r = 1$ and $k \approx 0,63$ align closely with the optimal values derived in Section~2. The measured external surface of $585,7 \,\mathrm{m}^2$ is virtually identical to the calculated optimum, with differences below $0,1 \,\mathrm{m}^2$ across all scenarios. Consequently, the compactness measure is equal to unity, confirming that the design is already optimal in practice. This case demonstrates that when geometric proportions are carefully balanced, near-perfect compactness can be achieved without compromise.

Comparing the three buildings highlights the practical lessons of the optimization framework. Elongated footprints, as in House A, tend to reduce efficiency by increasing surface area relative to volume. Extremely low heights, as in House B, similarly produce unfavorable envelope-to-volume ratios that result in cost penalties. Conversely, House C illustrates that compact, square-based designs with balanced slenderness achieve the theoretical optimum, combining structural economy with energy efficiency. These results confirm that compactness optimization is not merely a theoretical construct but a practical design guideline. Small adjustments to proportions in Houses A and B could have led to substantial savings in construction and operating costs, whereas House C exemplifies how theory and practice can converge in a well-proportioned architectural form.

\section{Software application}

To support the practical use of the theoretical results presented in this work, we have developed a free and openly accessible computer application dedicated to the optimization of hip roof houses.  
The program has been designed as an intuitive and user-friendly tool that enables architects, engineers, and students to directly apply the optimization procedures described in this study to real design tasks.

\begin{figure}[t]
\centering
\includegraphics[width=16 cm]{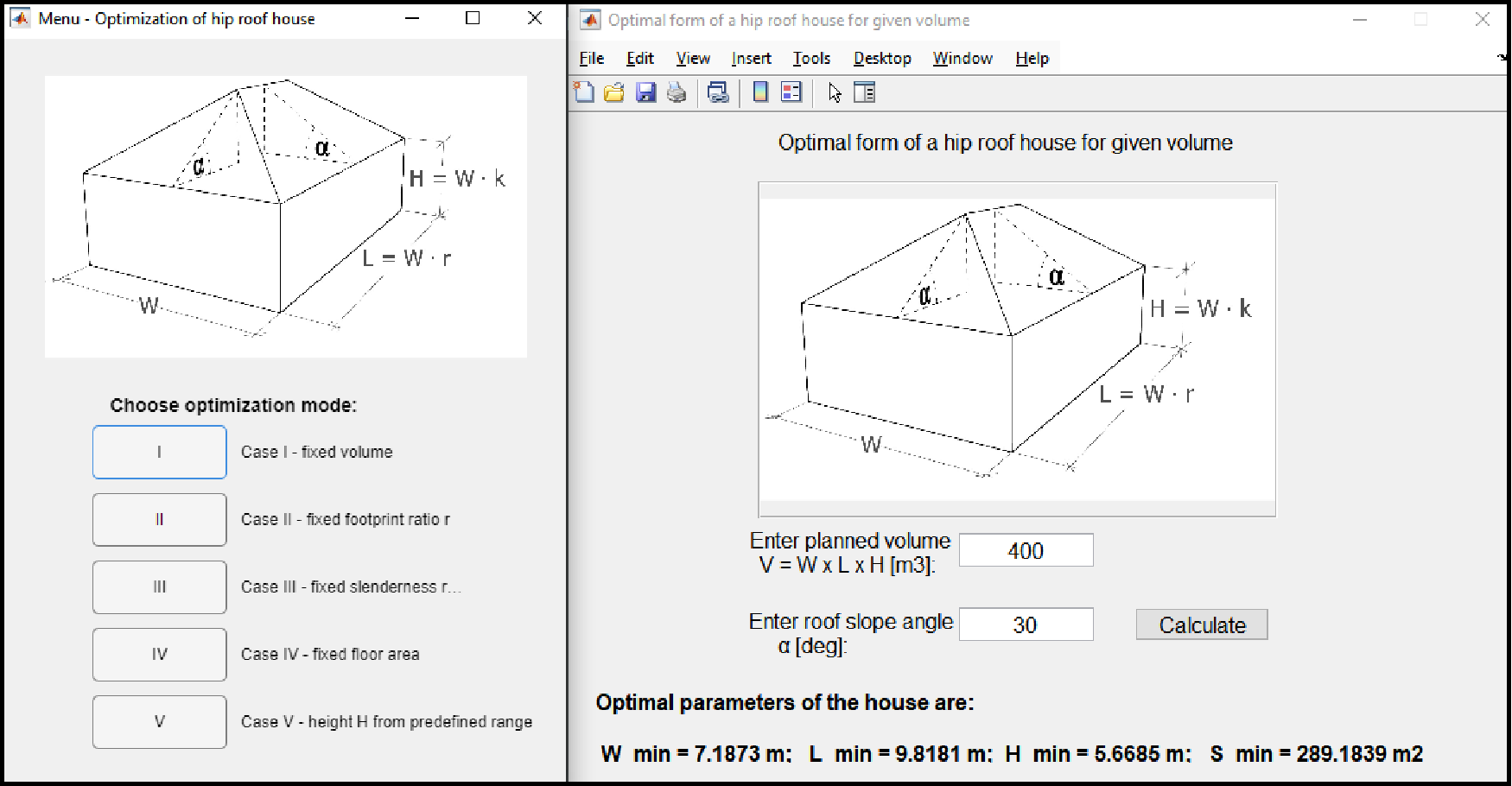}
\caption{Screenshot of the developed and shared computer program designed to perform optimization analyses of hip roof houses. }
\label{programy}
\end{figure}

Upon launching the application, the user is presented with a start-up window that allows the selection of the optimization mode (see Fig.~\ref{programy}, left panel).  
Five optimization scenarios analyzed in this paper are available. 
After choosing the desired mode, the program opens the corresponding analysis window, where the user can specify the input parameters (such as planned volume or roof slope angle) and obtain the resulting optimal house dimensions, $(W_{\min}, L_{\min}, H_{\min})$, together with the minimal external surface $S_{\min}$ (see Fig.~\ref{programy}, right panel).

The application is available free of charge at the Rokita-Projekt Architectural Office website:  
\url{https://rokita-projekt.pl/hip_roof_house_optimization.zip}.  

This software bridges the gap between mathematical optimization and architectural practice.  
It can be employed to explore alternative design solutions, verify the compactness and efficiency of proposed house forms, and assess potential cost savings.  
Beyond professional applications, it also serves as an educational tool, demonstrating how optimization methods can be directly translated into practical design guidelines.

The program was written in Matlab R2024b environment. MATLAB Runtime (R2024b) package is required to run the program. The package installer is freely available on the official MathWorks website.

\section{Conclusion}

The research presented in this paper has developed a rigorous mathematical framework for the optimization of hip roof houses, linking geometric proportions with sustainability-driven performance criteria. By systematically analyzing several scenarios — fixed volume, fixed footprint ratio, fixed slenderness, fixed floor area, and height constrained within a feasible range — we established explicit formulas for the optimal dimensions of hip roof houses. These formulas not only minimize the external surface of the building envelope but also provide clear guidelines for achieving material efficiency, cost reduction, and improved energy performance. The derived relationships confirmed that square-based footprints combined with balanced slenderness ratios consistently yield the most compact and efficient forms, while deviations toward elongated or flattened proportions significantly increase external surface area, and thereby both construction costs and energy demands. The case studies demonstrated that compactness optimization is not merely a theoretical construct but a practical guideline for architects and engineers, where minor design corrections have the potential to deliver substantial long-term sustainability gains.

The importance of this research lies in its dual contribution. On the one hand, it provides a mathematically rigorous and scale-independent compactness measure that serves as an objective criterion for evaluating and comparing architectural forms. On the other hand, it bridges geometry with sustainability by directly relating optimization of shape to reductions in resource consumption and energy use. In an era of increasing emphasis on sustainable construction, energy efficiency, and cost-effectiveness, the methods and results of this study enrich architectural design with a robust analytical tool. Furthermore, the accompanying software application developed as part of this research ensures that the theoretical outcomes can be seamlessly applied in practice. By making optimization accessible to architects, engineers, and students, the tool transforms abstract results into a decision-making aid that enhances both professional practice and education.

Despite the robustness of the presented results, this research opens multiple possibilities for future development. First, the optimization framework could be extended to more complex roof geometries, including gable-hip hybrids, dormer extensions, or multi-wing layouts, which are common in residential design. Second, while the present work focused on minimizing the external surface as a proxy for material and energy efficiency, future models could integrate additional performance indicators, such as solar gain optimization, daylight access, natural ventilation potential, or embodied carbon metrics. Third, dynamic constraints such as seismic performance, wind resistance, and structural load-bearing capacity could be explicitly incorporated into the optimization, thereby providing even more comprehensive and realistic design guidelines. Another promising direction involves the application of multi-objective optimization, where trade-offs between compactness, aesthetic preferences, functional requirements, and sustainability targets are systematically balanced. 

In conclusion, this study demonstrates that the intersection of geometry, mathematics, and sustainability creates powerful opportunities for improving architectural design. The systematic optimization of hip roof houses provides both theoretical insight and practical tools that are directly applicable to real construction projects. By embedding mathematical rigor into architectural practice, we can achieve forms that are simultaneously elegant, structurally efficient, economically viable, and environmentally responsible.


\section*{Data accessibility}
This article has no additional data.

\section*{Authors’ contributions}
All authors contributed equally to the conception of the research, as well as to the interpretation of results. E.R-M. designed the research, contributed to the methodological development of the research and wrote the original draft. B.G. reviewed the manuscript and provided final approval for publication. M.M. performed formal analysis and wrote the original draft.

\section*{Competing interests} 
The authors declare no competing interests.

\section*{Funding} 
This article has no funding.

\section*{Acknowledgements}
This article has no acknowledgements.

\bibliographystyle{plain}
\bibliography{bibliography}

@article{d2019compactness,
  title={A compactness measure of sustainable building forms},
  author={D’Amico, Bernardino and Pomponi, Francesco},
  journal={Royal Society Open Science},
  volume={6},
  number={6},
  pages={181265},
  year={2019},
  publisher={The Royal Society}
}

@misc{hip_roof_www,
  author       = {},
  title        = {},
  howpublished = {https://blog.newhomesource.com/hip-roofs},
  note         = {},
  year         = {2024},
  key          = {hip_roof_www}
}

@article{depecker2001design,
  title={Design of buildings shape and energetic consumption},
  author={Depecker, Patrick and Menezo, Christophe and Virgone, Joseph and Lepers, Stephane},
  journal={Building and Environment},
  volume={36},
  number={5},
  pages={627--635},
  year={2001},
  publisher={Elsevier}
}

@book{hegger2012energy,
  title={Energy manual: sustainable architecture},
  author={Hegger, Manfred and Fuchs, Matthias and Stark, Thomas and Zeumer, Martin},
  year={2012},
  publisher={Walter de Gruyter}
}

@book{Steadman2000,
  title={The Evolution of Designs: Biological Analogy in Architecture and the Applied Arts},
  author={Steadman, Philip},
  publisher={Routledge},
  year={2000}
}

@book{Steadman2014,
  title={Building Types and Built Forms},
  author={Steadman, Philip},
  publisher={Troubador Publishing},
  year={2014}
}

@article{Julia2017,
  title={Validation of a Bayesian-based method for defining residential archetypes in urban building energy models},
  author={Sokol, Julia and Davila, Carlos Cerezo and Reinhart, Christoph F},
  journal={Energy and Buildings},
  volume={134},
  pages={11--24},
  year={2017},
  publisher={Elsevier}
}

@article{Hargreaves2017,
  title={Representing the dwelling stock as 3D generic tiles estimated from average residential density},
  author={Hargreaves, AJ},
  journal={Computers, Environment and Urban Systems},
  volume={54},
  pages={280--300},
  year={2015},
  publisher={Elsevier}
}

@book{Angel2012,
  title={Planet of Cities},
  author={Angel, Shlomo},
  publisher={Lincoln Institute of Land Policy},
  year={2012}
}

@techreport{UN2014,
  title={World Urbanization Prospects: 2014 Revision},
  author={ },
  year={2014},
  institution={Department of Economic and Social Affairs, United Nations}
}

@misc{LondonStrategy2018,
  title={London Environment Strategy},
  author={{Greater London Authority}},
  year={2018},
  howpublished={https://www.london.gov.uk/what-we-do/environment/london-environment-strategy}
}

@misc{California2017,
  title={Buy Clean California Act},
  author={{State of California}},
  year={2017},
  note={Assembly Bill No. 262}
}

@article{Ramesh2010,
  title={Life cycle energy analysis of buildings: An overview},
  author={Ramesh, T. and Prakash, Ravi and Shukla, K.K.},
  journal={Energy and Buildings},
  volume={42},
  number={10},
  pages={1592--1600},
  year={2010}
}

@article{Martin1972,
  title={Architect’s approach to architecture},
  author={Martin, Leslie},
  journal={RIBA},
volume={74},
  pages={191--200},
  year={1967}
}

@article{March1976,
  title={The geometry of environment},
  author={March, Lionel and Steadman, Philip},
  year={1971},
  publisher={Routledge}
}

@article{Steemers2003,
  title={Energy and the city: density, buildings and transport.},
  author={Steemers, Koen},
  journal={Energy and Buildings},
  volume={35},
  pages={3--14},
  year={2003}
}

@article{Catalina2011,
  title={Study on the impact of the building form on the energy consumption},
  author={Catalina, Tiberiu and Virgone, Joseph and Iordache, Vlad},
  booktitle={Building Simulation 2011},
  volume={12},
  pages={1726--1729},
  year={2011},
  organization={IBPSA}
}

@article{Ourghi2007,
  title={Simplified analysis of the effect of shape on annual energy use},
  author={Ourghi, R. and Al-Anzi, A. and Krarti, M.},
  journal={Energy Conversion and Management},
  volume={48},
  pages={300--305},
  year={2007}
}

@article{Schlueter2009,
  title={Building information model based energy/exergy performance assessment in early design stages},
  author={Schlueter, Arno and Thesseling, Frank},
  journal={Automation in construction},
  volume={18},
  number={2},
  pages={153--163},
  year={2009},
  publisher={Elsevier}
}

@article{Jedrzejuk2000,
  title={Optimization of shape and functional structure of buildings as well as heat source utilization. Basic theory},
  author={Jedrzejuk, Hanna and Marks, Wojciech},
  journal={Building and Environment},
  volume={37},
  number={12},
  pages={1379--1383},
  year={2002},
  publisher={Elsevier}
}

@article{Jedrzejuk2002,
  title={Optimization of shape and functional structure of buildings as well as heat source utilisation. Partial problems solution},
  author={Jedrzejuk, Hanna and Marks, Wojciech},
  journal={Building and Environment},
  volume={37},
  number={11},
  pages={1037--1043},
  year={2002},
  publisher={Elsevier}
}

@article{Hachem2012,
  title={Parametric investigation of geometric form effects on solar potential of housing units},
  author={Hachem, Caroline and Athienitis, Andreas and Fazio, Paul},
  journal={Solar Energy},
  volume={85},
  number={9},
  pages={1864--1877},
  year={2011},
  publisher={Elsevier}
}

@article{Hachem2016,
  title={Investigation of solar potential of housing units in different neighborhood designs},
  author={Hachem, Caroline and Athienitis, Andreas and Fazio, Paul},
  journal={Energy and Buildings},
  volume={43},
  number={9},
  pages={2262--2273},
  year={2011},
  publisher={Elsevier}
}

@article{Okeil2010,
  title={A holistic approach to energy efficient building forms},
  author={Okeil, Ahmad},
  journal={Energy and buildings},
  volume={42},
  number={9},
  pages={1437--1444},
  year={2010},
  publisher={Elsevier}
}

@article{Caruso2013,
  title={Optimal theoretical building form to minimize direct solar irradiation},
  author={Caruso, Gianpiero and Fantozzi, Fabio and Leccese, Francesco},
  journal={Solar Energy},
  volume={97},
  pages={128--137},
  year={2013},
  publisher={Elsevier}
}

@article{Jin2017,
  title={Optimization of a free-form building shape to minimize external thermal load using genetic algorithm},
  author={Jin, Jeong-Tak and Jeong, Jae-Weon},
  journal={Energy and Buildings},
  volume={85},
  pages={473--482},
  year={2014},
  publisher={Elsevier}
}

@article{Vartholomaios2017,
  title={A parametric sensitivity analysis of the influence of urban form on domestic energy consumption for heating and cooling in a Mediterranean city},
  author={Vartholomaios, Aristotelis},
  journal={Sustainable cities and society},
  volume={28},
  pages={135--145},
  year={2017},
  publisher={Elsevier}
}

@book{boyd2004convex,
  title={Convex optimization},
  author={Boyd, Stephen P and Vandenberghe, Lieven},
  year={2004},
  publisher={Cambridge university press}
}

@article{stiny1980,
  title={Introduction to shape and shape grammars},
  author={Stiny, George},
  journal={Environment and planning B: planning and design},
  volume={7},
  number={3},
  pages={343--351},
  year={1980},
  publisher={SAGE Publications Sage UK: London, England}
}

@book{mitchell1990,
  author    = {Mitchell, William J.},
  title     = {The Logic of Architecture: Design, Computation, and Cognition},
  year      = {1990},
  publisher = {MIT Press},
  address   = {Cambridge, MA},
  isbn      = {978-0262132555}
}

@book{knight1994,
  author    = {Knight, Terry},
  title     = {Transformations in Design: A Formal Approach to Stylistic Change and Innovation in the Visual Arts},
  year      = {1994},
  publisher = {Cambridge University Press},
  address   = {Cambridge, UK},
  isbn      = {978-0521421836}
}

@article{oxman2006theory,
  author    = {Oxman, Rivka},
  title     = {Theory and design in the first digital age},
  journal   = {Design Studies},
  volume    = {27},
  number    = {3},
  pages     = {229--265},
  year      = {2006},
  publisher = {Elsevier},
  doi       = {10.1016/j.destud.2005.11.002},
  issn      = {0142-694X}
}

\end{document}